\documentclass[11pt]{article}
\usepackage{amssymb}
\usepackage{amscd}
\usepackage{amsmath}
\usepackage{comment}
\input amssym.def
\input amssym

\newtheorem{theorem}{Theorem}[section]
\newtheorem{corollary}[theorem]{Corollary}
\newtheorem{lemma}[theorem]{Lemma}

\newtheorem{conjecture}[theorem]{Conjecture}

\newtheorem{remark}[theorem]{Remark}
\newtheorem{definition}[theorem]{Definition}
\newtheorem{proposition}[theorem]{Proposition}

\usepackage{amssymb} 
\title{Principal
subspaces of higher-level standard
$\widehat{\mathfrak{sl}(n)}$-modules} \author{ Christopher Sadowski}
\date{}
\begin{document}
\maketitle

\renewcommand{\theequation}{\thesection.\arabic{equation}}
\renewcommand{\thetheorem}{\thesection.\arabic{theorem}}
\setcounter{equation}{0} \setcounter{theorem}{0}
\setcounter{section}{0}
\begin{abstract}
Using completions of certain universal enveloping algebras, we provide
a natural setting for families of defining relations
for the principal subspaces of standard modules for
untwisted affine Lie algebras. We also use the theory of vertex operator algebras
and intertwining operators to construct exact sequences among principal subspaces of certain
standard $\widehat{\mathfrak{sl}(n)}$-modules, $n \ge 3$. As 
a consequence, we obtain the multigraded dimensions of the principal subspaces
$W(k_1\Lambda_1 + k_2 \Lambda_2)$ and $W(k_{n-2}\Lambda_{n-2} + k_{n-1} \Lambda_{n-1})$.
This generalizes earlier
work by Calinescu on principal subspaces of 
standard $\widehat{\mathfrak{sl}(3)}$-modules.
\end{abstract}

\section{Introduction}

Principal subspaces of standard modules for affine Lie algebras were
introduced and studied by Feigin and
Stoyanovsky in \cite{FS1}--\cite{FS2}.
In their work, motivated by the earlier work by Lepowsky and
Primc \cite{LP}, Feigin and Stoyanovsky discovered that the
multigraded dimensions
(generating functions of dimensions of homogeneous subspaces)
of principal subspaces of standard
$\widehat{\mathfrak{sl}(2)}$-modules are related to the
Rogers-Ramanujan partition identities, and more generally, the
Andrews-Gordon identities (cf. \cite{A}). Multigraded dimensions for a
more general class of principal subspaces were later studied by
Georgiev in \cite{G}, where combinatorial bases were constructed for
the principal subspaces of certain standard
$\widehat{\mathfrak{sl}(n+1)}$-modules, $n \ge 1$ (written for
brevity as $\widehat{\mathfrak{sl}(n)}$ in the title).
More recently, combinatorial
bases have been constructed for principal subspaces in more general
lattice cases (\cite{P}, \cite{MiP}), for the principal subspaces of the
vacuum standard modules for the affine Lie algebra $B_2^{(1)}$
\cite{Bu}, for principal subspaces in the quantum
$\widehat{\mathfrak{sl}(n+1)}$-case \cite{Ko}, and for certain 
substructures of principal subspaces (\cite{Pr},
\cite{J1}--\cite{J3}, \cite{T1}--\cite{T4}, \cite{Ba}, \cite{JPr}).

In \cite{CLM1}--\cite{CLM2}, Capparelli, Lepowsky, and Milas
interpreted the Rogers-Ramanujan and Rogers-Selberg recursions in terms
of the multigraded dimensions of the principal subspaces of the standard
$\widehat{\mathfrak{sl}(2)}$-modules by using the
vertex-algebraic structure of these
modules, along with intertwining
operators among them, to construct exact sequences.
In \cite{CLM1}--\cite{CLM2} (as in
\cite{FS1}--\cite{FS2}), the authors assumed certain presentations
(generators and defining relations) for the principal subspaces of the
standard $\widehat{\mathfrak{sl}(2)}$-modules,
presentations that can be derived from \cite{LP}; the nontrivial
part is the completeness of the relations. The question of
proving in an a priori way that the relations assumed in these works were indeed a
complete set of defining relations for the principal subspaces of the
standard $\widehat{\mathfrak{sl}(2)}$-modules was later addressed by
Calinescu, Lepowsky, and Milas in \cite{CalLM1}--\cite{CalLM2}, where
the authors gave such an a priori proof.
These results were extended to the level $1$ standard
$\widehat{\mathfrak{sl}(n+1)}$-modules by Calinescu in \cite{C2}, and
later to the level $1$ standard modules for the
untwisted affine Lie algebras of type $ADE$ in \cite{CalLM3}.
In both \cite{C2} and \cite{CalLM3}, the
authors proved that the multigraded dimension of the principal
subspace of the vacuum module satisfies a certain recursion, and using
this recursion they found the multigraded dimensions of the principal
subspaces of all the level $1$ standard modules.
 
In the work \cite{C1}, Calinescu considered the principal subspaces of
certain higher level standard $\widehat{\mathfrak{sl}(3)}$-modules. In
this work, she conjecturally assumed presentations for certain principal subspaces,
and using the theory of vertex operator algebras and intertwining
operators, she constructed exact sequences among these principal
subspaces. Using these exact sequences, along with the multigraded
dimensions in \cite{G}, Calinescu was able to find the multigraded
dimensions of principal subspaces which had not previously been
studied. In \cite{S} the presentations for the principal subspaces of
all the standard $\widehat{\mathfrak{sl}(3)}$-modules were proved
(including those assumed in \cite{C1}).

Our main result in the present
work is a natural generalization of \cite{C1} to the case of
$\widehat{\mathfrak{sl}(n+1)}$, $n\ge 2$.
Although our methods recover the same information
as in \cite{CLM1}--\cite{CLM2} when $n =1$, we take $n \ge 2$ 
for notational convenience.
In the case where $n=2$, we recover the results in \cite{C1} 
with a slight variant of the methods.
Using the work of \cite{G}, we 
provide exact sequences among
principal subspaces of certain standard
$\widehat{\mathfrak{sl}(n+1)}$-modules. As a consequence, we
obtain previously unknown multigraded dimensions of principal subspaces.
In addition, as in \cite{C1}, we conjecturally assume presentations
for certain principal subspaces, and use these to obtain 
exact sequences among a more general class of
principal subspaces of certain standard
$\widehat{\mathfrak{sl}(n+1)}$-modules.
To state our main result, we let
$\Lambda_0,\dots,\Lambda_n$ denote the fundamental weights of
$\widehat{\mathfrak{sl}(n+1)}$. The dominant integral weights $\Lambda$
of $\widehat{\mathfrak{sl}(n+1)}$ are $k_0\Lambda_0 + \dots
+ k_n \Lambda_n$ for $k_0,\dots,k_n \in \mathbb{N}$, and we use
$L(\Lambda)$ to denote the standard module with highest weight
$\Lambda$, $W(\Lambda)$ to denote its principal subspace, and
$\chi'_{W(\Lambda)}(x_1,\dots,x_n,q)$ to denote its multigraded
dimension. Our result states:

\begin{theorem}\label{introseq}
Let $k \ge 1$.
For $k_1, k_{2}, k_{n-1}, k_n \in \mathbb{N}$
such that $k_1 + k_{2} = k_{n-1}+ k_n = k$ and $k_1 > 0$ and $k_n > 0$ the sequences
\begin{eqnarray} \label{introseq1}
\lefteqn{0 \longrightarrow W(k_1 \Lambda_1 + k_2 \Lambda_2)
  \stackrel{e^{\otimes k}_{\omega_1}} \longrightarrow} \\ && W(k_1
\Lambda_0 + k_2 \Lambda_1) \stackrel{1^{\otimes{k_1-1}} \otimes{\cal
    Y}_{c}(e^{\lambda_1},x)\otimes 1^{\otimes k_{2}}}
\longrightarrow \nonumber \\ && \hspace{2em} W((k_1-1)
\Lambda_0+(k_2+1) \Lambda_1) \longrightarrow 0 \nonumber
\end{eqnarray}
and
\begin{eqnarray}\label{introseq2}
\lefteqn{0 \longrightarrow W(k_{n-1} \Lambda_{n -1} + k_{n}
  \Lambda_{n}) \stackrel{e^{\otimes k}_{\omega_n}} \longrightarrow}
\\ && W(k_n \Lambda_0 + k_{n-1} \Lambda_n)
\stackrel{1^{\otimes{k_i-1}} \otimes{\cal
    Y}_{c}(e^{\lambda_i},x)\otimes 1^{\otimes k_{i+1}}}
\longrightarrow \nonumber \\ && \hspace{2em} W((k_n-1)
\Lambda_0+(k_{n-1}+1) \Lambda_n) \longrightarrow 0 \nonumber
\end{eqnarray}
are exact.
\end{theorem}

More generally, conjecturally assuming certain presentations, we obtain:

\begin{theorem} 
Let $k \ge 1$. 
For any $i$ with $1 \leq i \leq n-1$ and $k_i, k_{i+1} \in \mathbb{N}$ such that 
$k_i + k_{i+1} = k$, the sequences
\begin{eqnarray} 
\lefteqn{W(k_i \Lambda_i + k_{i+1} \Lambda_{i+1}) \stackrel{\phi_i}
  \longrightarrow} \\ && W(k_i \Lambda_0 +k_{i+1} \Lambda_i)
\stackrel{1^{\otimes{k_i-1}} \otimes{\cal
    Y}_{c}(e^{\lambda_i},x)\otimes 1^{\otimes k_{i+1}}}
\longrightarrow \nonumber \\ && \hspace{2em} W((k_i-1)
\Lambda_0+(k_{i+1}+1) \Lambda_i) \longrightarrow 0 \nonumber
\end{eqnarray}
when $k_i \ge 1$, and
\begin{eqnarray} 
\lefteqn{W(k_i \Lambda_i + k_{i+1} \Lambda_{i+1}) \stackrel{\psi_i}
  \longrightarrow} \\ && W(k_{i+1} \Lambda_0 + k_i \Lambda_{i+1})
\stackrel{1^{\otimes{k_{i+1}-1}} \otimes{\cal
    Y}_{c}(e^{\lambda_{i+1}},x)\otimes 1^{\otimes k_{i}}}
\longrightarrow \nonumber \\ && \hspace{2em} W((k_{i+1}-1)
\Lambda_0+(k_i+1) \Lambda_{i+1}) \longrightarrow 0 \nonumber
\end{eqnarray}
when $k_{i+1} \ge 1$, are exact.
\end{theorem}

The maps $\phi_i$, $\psi_i$, $e_{\omega_1}^{\otimes k_1}$,
$e_{\omega_n}^{\otimes k_n}$, and ${\cal Y}_{c}(e^{\lambda_i},x)$ are
maps naturally arising from the lattice construction of the
level $1$ standard
modules and intertwining operators among these modules.  As a
consequence Theorem \ref{introseq}, we obtain results about multigraded
dimensions, and we have the following theorem and its corollary:
\begin{theorem} \label{introtheorem}
Let $k \ge 1$.
Let $k_1, k_2, k_{n-1}, k_n \in \mathbb{N}$ with $k_1 \ge 1$ and $k_n \ge 1$,
such that $k_1 + k_2 = k$ and $k_{n-1} + k_n = k$. Then
\begin{eqnarray} 
\lefteqn{\chi'_{W(k_1\Lambda_1 + k_2\Lambda_2)}(x_1,\dots,x_n,q) = \nonumber }\\
&=&x_1^{-k_1}\chi'_{W((k_1-1)\Lambda_0 + (k_2+1)\Lambda_1)}(x_1q^{-1},x_2q, x_3\dots,x_n,q)\\
&& - 
x_1^{-k_1}\chi'_{W(k_1\Lambda_0 + k_2\Lambda_1)}(x_1q^{-1},x_2q,x_3,\dots,x_n,q)\nonumber
\end{eqnarray}
and
\begin{eqnarray}
\lefteqn{\chi'_{W(k_{n-1}\Lambda_{n-1} + k_n \Lambda_n)}(x_1,\dots,x_n,q) = \nonumber} \\
&=&x_n^{-k_n}\chi'_{W((k_n-1)\Lambda_0 + (k_{n-1}+1)\Lambda_n)}(x_1,\dots,x_{n-1}q,x_nq^{-1},q)\\
&& - 
x_n^{-k_n}\chi'_{W(k_n\Lambda_0 + k_{n-1}\Lambda_n)}(x_1,\dots,x_{n-1}q,x_nq^{-1},q).\nonumber
\end{eqnarray}
\end{theorem}
Theorem \ref{introtheorem} immediately gives us:
\begin{corollary}\label{introgdim}
In the setting of Theorem \ref{introtheorem}, we have that
$$
\chi'_{W(k_1\Lambda_1 + k_2\Lambda_2)}(x_1,\dots,x_n,q)=
$$
$$
=\sum_{}\;
\bigg(\frac{q^{{r_{1}^{(1)}}^{2}+ \ldots + {r_{1}^{(k)}}^{2} +
\sum_{t=k_1+1}^{k} r_{1}^{(t)} + \sum_{t=1}^k r_2^{(t)} - r_1^{(t)} }
(1-q^{r_1^{(k_1)}})}{(q)_{r_{1}^{(1)} -
r_{1}^{(2)}} \ldots (q)_{r_{1}^{(k - 1)} - r_{1}^{(k)}}
(q)_{r_{1}^{(k)}}}\bigg)
\bigg( \frac{q^{{r_{2}^{(1)}}^{2} + \ldots + {r_{2}^{(k)}}^{2} -
r_{2}^{(1)} r_{1}^{(1)} - \ldots - r_{2}^{(k)}r_{1}^{(k)}
 }}{(q)_{r_{2}^{(1)}
- r_{2}^{(2)}} \ldots (q)_{r_{2}^{(k - 1)}
- r_{2}^{(k)}} (q)_{r_{2}^{(k)}}}\bigg) \times$$
$$\times\cdots\times \bigg(
\frac{q^{{r_{n}^{(1)}}^{2} + \ldots + {r_{n}^{(k)}}^{2} -
r_{n}^{(1)} r_{n - 1}^{(1)} - \ldots - r_{n}^{(k)}r_{n - 1}^{(k)}
}}{(q)_{r_{n}^{(1)}
- r_{n}^{(2)}} \ldots (q)_{r_{n}^{(k - 1)}
- r_{n}^{(k)}} (q)_{r_{n}^{(k)}}}\bigg) x_1^{-k_1 + \sum_{i=1}^k r_1^{(i)}}\cdots x_n^{\sum_{i=1}^n r_n^{(i)}}$$
and
$$
\chi'_{W(k_{n-1}\Lambda_{n-1} + k_{n}\Lambda_n)}(x_1,\dots,x_n,q)=
$$
$$
=\sum_{}\;
\bigg(\frac{q^{{r_{1}^{(1)}}^{2}+ \ldots + {r_{1}^{(k)}}^{2} 
}}{(q)_{r_{1}^{(1)} -
r_{1}^{(2)}} \ldots (q)_{r_{1}^{(k - 1)} - r_{1}^{(k)}}
(q)_{r_{1}^{(k)}}}\bigg)
\bigg(\frac{q^{{r_{2}^{(1)}}^{2} + \ldots + {r_{2}^{(k)}}^{2} -
r_{2}^{(1)} r_{1}^{(1)} - \ldots - r_{2}^{(k)}r_{1}^{(k)}+
}}{(q)_{r_{2}^{(1)}
- r_{2}^{(2)}} \ldots (q)_{r_{2}^{(k - 1)}
- r_{2}^{(k)}} (q)_{r_{2}^{(k)}}}\bigg)\times$$
$$\times\cdots\times
\bigg(\frac{q^{{r_{n}^{(1)}}^{2} + \ldots + {r_{n}^{(k)}}^{2} -
r_{n}^{(1)} r_{n - 1}^{(1)} - \ldots - r_{n}^{(k)}r_{n - 1}^{(k)}+
\sum_{t=k_n+1}^{k} r_{n}^{(t)}
}
}{(q)_{r_{n}^{(1)}
- r_{n}^{(2)}} \ldots (q)_{r_{n}^{(k - 1)}
- r_{n}^{(k)}} (q)_{r_{n}^{(k)}}}\bigg) \times
$$
$$\times q^{\sum_{t=1}^kr_{n-1}^{(t)}-r_{n}^{(t)}}(1-q^{r_n^{(k_n)}})
x_1^{\sum_{i=1}^k r_1^{(i)}}\cdots x_n^{-k_n + \sum_{i=1}^n r_n^{(i)}}$$
where the sums are taken over decreasing sequences 
$r_j^{(1)} \ge r_j^{(2)} \ge \dots \ge r_j^{(k)} \ge 0$ for each $j=1,\dots,n$.
\end{corollary}

The expressions in Corollary \ref{introgdim} can also be written as follows:
As in \cite{G}, for $s= 1,\dots,k-1$ and $i=1,\dots,n$, set $p_i^{(s)} = r_i^{(s)} - r_i^{(s+1)}$,
and set $p_i^{(k)} = r_i^{(k)}$. Also, let $(A_{lm})_{l,m = 1}^n $ be the Cartan matrix of
$\mathfrak{sl}(n+1)$ and $B^{st} := \mathrm{min}\{s,t\}$, $1 \le s,t \le k$. Then,
$$\chi'_{W(k_1\Lambda_1 + k_2\Lambda_2)}(x_1,\dots,x_n,q)=$$
$$\sum_{\begin{array}{c}{\scriptstyle
p_{1}^{(1)} ,\ldots ,  p_{1}^{(k)} \geq 0}\\{\scriptstyle 
\vdots}\\{\scriptstyle p_{n}^{(1)}, \ldots
, p_{n}^{(k)} \geq 0 }\end{array}}
\frac{q^{\frac{1}{2}\sum_{l,m = 1}^n\sum_{s,t=1}^lA_{lm}B^{st}p_l^{(s)} p_m^{(t)}}}{
\prod_{i=1}^n \prod_{s=1}^k (q)_{{p_i}^{(s)}}}q^{\widetilde{p_1}}
q^{\sum_{t=1}^k p_2^{(t)} + \dots + p_2^{(k)} - p_1^{(t)} - \dots - p_1^{(k)}}\times$$
$$\times 
(1-q^{p_1^{(k_1)} + \dots + p_1^{(k)}})x_1^{-k_1}\prod_{i=1}^n x_i^{\sum_{s=1}^k sp_i^{(s)}}
$$
where $\widetilde{p_1} = p_1^{(k_1+1)} + 2p_1^{(k_1+2)} + \dots + k_2p_1^{(k)}$ and 
$$\chi'_{W(k_{n-1}\Lambda_{n-1} + k_n\Lambda_n)}(x_1,\dots,x_n,q)=$$
$$\sum_{\begin{array}{c}{\scriptstyle
p_{1}^{(1)} ,\ldots ,  p_{1}^{(k)} \geq 0}\\{\scriptstyle 
\vdots}\\{\scriptstyle p_{n}^{(1)}, \ldots
, p_{n}^{(k)} \geq 0 }\end{array}}
\frac{q^{\frac{1}{2}\sum_{l,m = 1}^n\sum_{s,t=1}^lA_{lm}B^{st}p_l^{(s)} p_m^{(t)}}}{
\prod_{i=1}^n \prod_{s=1}^k (q)_{{p_i}^{(s)}}}q^{\widetilde{p_n}}
q^{\sum_{t=1}^k p_{n-1}^{(t)} + \dots + p_{n-1}^{(k)} - p_n^{(t)} - \dots - p_n^{(k)}}
\times$$ $$
\times
(1-q^{p_n^{(k_n)} + \dots + p_n^{(k)}}) x_n^{-k_n}\prod_{i=1}^n x_i^{\sum_{s=1}^k sp_i^{(s)}}
$$
where $\widetilde{p_n} = p_n^{(k_n+1)} + 2p_n^{(k_n+2)} + \dots + k_{n-1}p_1^{(k)}$.

Similar multigraded dimensions for different variants of principal
subspaces have been studied in \cite{AKS} and \cite{FFJMM}.
Modularity
properties of certain multigraded dimensions, in the context of
principal subspaces of standard modules, have been studied in \cite{St},
\cite{WZ}, and more recently in
\cite{BCFK}.  In \cite{CalLM4}, the authors have initiated the study
of principal subspaces for standard modules for twisted affine Lie
algebras, extending the past work of \cite{CLM1}--\cite{CLM2},
\cite{CalLM1}--\cite{CalLM3} to the case of the level $1$ standard module
for the twisted affine Lie algebra $A_2^{(2)}$.

In \cite{CLM1}--\cite{CLM2}, \cite{CalLM1}--\cite{CalLM3}, and \cite{S},
the annihilator of the highest weight vector of each principal subspace
is written in terms of certain elements of 
$U(\widehat{\mathfrak{g}})$ which, when viewed as operators, annihilate the highest weight vector.
An important set of these operators arises from certain null vector identities
given by powers of vertex operators and are written as infinite
formal sums of elements of $U(\widehat{\mathfrak{g}})$ also viewed as operators.
The ideals which annihilate the highest weight vectors 
can be expressed using operators defined by certain truncations of these formal sums, in 
order to view these operators as
elements of $U(\widehat{\mathfrak{g}})$.

The remainder of the present work, including the Appendix,
focuses on providing a more natural (without such truncations)
setting for the annihilating ideals in order to give presentations of the
principal subspaces of the standard modules. A construction of a
completion of a universal enveloping algebra is used to give
natural presentations for the defining annihilating ideals of
principal subspaces.  This completion was discussed in
\cite{C1}--\cite{C2} and \cite{CalLM3}, but the details
of this construction were not supplied.
 We also prove various properties of
this completion and the defining ideals for principal subspaces,
including their more natural 
definition inside this completion. These completions may be
generalized to the twisted setting used in \cite{CalLM4}
(as in \cite{LW}, where similar completions were originally constructed
in a general twisted or untwisted setting).

We now give a brief outline of the structure of the present work. In
Section 2, we recall certain vertex operator constructions of the
standard $\widehat{\mathfrak{sl}(n+1)}$-modules and intertwining
operators among these modules.  In Section 3, we recall the definition
of the principal subspace of a standard module, along with the
definition of multigraded dimension. We also recall the
presentations of the principal subspaces of standard
$\widehat{\mathfrak{sl}(n+1)}$-modules, including known
and conjectured presentations. In Section 4, we reformulate these
presentations in terms of a completion of a certain universal
enveloping algebra. In Section 5, we construct exact sequences among
principal subspaces of certain standard
$\widehat{\mathfrak{sl}(n+1)}$-modules and obtain the multigraded
dimensions of $W(k_1\Lambda_1 + k_2\Lambda_2)$ and
$W(k_{n-1}\Lambda_{n-1} + k_n\Lambda_n)$. In the Appendix,
we recall a certain construction from \cite{LW} and use it to
construct the aforementioned completion of the universal enveloping
algebra.

{\bf Acknowledgment} I am very grateful to James Lepowsky for many insightful
discussions related to this work, along with many suggestions for improvements to this paper.
I would also like
to thank Antun Milas for helpful comments about the multigraded dimensions found in this work.
The author was partially supported by NSF grant PHY-0901237 and an
Excellence Fellowship from the Rutgers University Graduate School.

\section{Preliminaries}
\setcounter{equation}{0}

In this section we recall the vertex operator construction of the
standard $\widehat{\mathfrak{sl}(n+1)}$-modules.
We refer the reader to \cite{CalLM3} and \cite{S} for the
full details of these constructions as they are used here.  Let
$\mathfrak{g} = \mathfrak{sl}(n+1)$, where $n \ge 2$
(the rank of $\mathfrak{g}$), and fix a Cartan subalgebra $\goth{h}$
of $\mathfrak{g}$, a set of roots $\Delta$, a set of simple roots
$\{\alpha_1,\dots ,\alpha_n\}$, and a set of positive roots
$\Delta_+$.  We use $\langle \cdot, \cdot \rangle$ to denote the Killing
form and rescale it so that $\langle \alpha, \alpha \rangle =2$ for
each $\alpha \in \Delta$.  Using this form, we identify $\goth{h}$
with $\goth{h}^*$. Let $\lambda_1, \dots, \lambda_n \in \goth{h}
\simeq \goth{h}^*$ denote the fundamental weights of
$\goth{sl}(n+1)$. Recall that $\langle \lambda_i,\alpha_j \rangle =
\delta_{ij}$ for each $i,j=1,\dots, n$.  Denote by $Q= \sum_{i=1}^n
\mathbb{Z}\alpha_i$ and $P=\sum_{i=1}^n \mathbb{Z}\lambda_i$ the root
lattice and weight lattice of $\goth{sl}(n+1)$, respectively.

For each root $\alpha \in \Delta$, fix a root vector $x_{\alpha} \in
\goth{sl}(n+1)$.  Define
\begin{eqnarray*} 
\goth{n}= \sum_{\alpha \in \Delta_+} \mathbb{C} x_\alpha,
\end{eqnarray*}
a nilpotent subalgebra of $\goth{sl}(n+1)$.

We have the corresponding untwisted affine Lie algebra given by
\begin{eqnarray*}
\widehat{{\goth{sl}(n+1)}}= {\goth{sl}(n+1)} \otimes \mathbb{C}[t, t^{-1}] \oplus \mathbb{C}c,
\end{eqnarray*}
where $c$ is a non-zero central element and 
\begin{eqnarray*}
[ x \otimes t^m, y \otimes t^p ] = [x, y] \otimes t^{m+p} + m\langle x, y \rangle \delta _{m+p, 0} c
\end{eqnarray*}
for any $x, y \in {\goth{sl}(n+1)}$ and $m, p \in \mathbb{Z}$.  Adjoining
the degree operator $d$

we obtain the affine Kac-Moody Lie algebra
$\widetilde{{\goth{sl}(n+1)}}=\widehat{{\goth{sl}(n+1)}} \oplus
\mathbb{C}d$ (cf. \cite{K}). 
We extend our form $\langle \cdot,\cdot \rangle$ to $\goth{h} \oplus \mathbb{C}c \oplus \mathbb{C}d $. 

Using this form,
we may identify $\goth{h} \oplus \mathbb{C}c \oplus \mathbb{C}d$ with
$(\goth{h} \oplus \mathbb{C}c \oplus \mathbb{C}d)^{*}$.  The simple
roots of $\widehat{{\goth{sl}(n+1)}}$ are $\alpha _{0}$, $\alpha
_{1},\dots , \alpha _{n}$ and the fundamental weights of
$\widehat{{\goth{sl}(n+1)}}$ are $\Lambda_{0}, \Lambda _{1},\dots,
\Lambda _{n} $, given by
\begin{eqnarray*}
\alpha_0= c-(\alpha_1+\alpha_2 + \dots + \alpha_n)
\end{eqnarray*}
and
\begin{eqnarray*}
 \Lambda_0= d, \; \Lambda_i=\Lambda_0+\lambda_i
\end{eqnarray*}
for each $i=1, \dots, n$. For each weight $\Lambda_i$, $i=0,\dots,n$, we denote by
$L(\Lambda_i)$ the level $1$ standard module with highest weight weight $\Lambda_i$,
and we denote its highest weight vector by $v_{\Lambda_i}$.

We now recall the lattice vertex operator construction of the level $1$ standard
$\widehat{\goth{sl}(n+1)}$-modules. $\widehat{\goth{sl}(n+1)}$ has two
important subalgebras:
\begin{eqnarray*}
\widehat{\goth{h}} = \goth{h} \otimes \mathbb{C}[t, t^{-1}] \oplus
\mathbb{C}c
\end{eqnarray*}
and the Heisenberg subalgebra
\begin{eqnarray*}
  \widehat {\goth{h}}_{\mathbb{Z}} = \coprod _{m \in \mathbb{ Z}
    \setminus \{0\}} \goth{h} \otimes t^m \oplus \mathbb{C}c
\end{eqnarray*}
(in the notation of \cite{FLM}, \cite{LL}).

We use $U(\cdot)$ to denote the universal enveloping algebra.
 The induced module
\begin{eqnarray*} 
M(1)= U( \widehat{\goth{h}}) \otimes _{U(\goth{h} \otimes \mathbb{C}[t] \oplus \mathbb{C}c)}\mathbb{C}
\end{eqnarray*}
has a natural $\widehat{\goth{h}}$-module structure, where
$\goth{h}\otimes \mathbb{C}[t]$ acts trivially and $c$ acts as
identity on the one-dimensional module $ \mathbb{C}.$

We use the notation $\mathbb{C}[S]$ to denote the group algebra of $S$
for any $S \subset P$, and define
$$ V_P = M(1) \otimes \mathbb{C}[P],$$
$$V_Q= M(1) \otimes \mathbb{C}[Q]$$ 
and we set
$$ V_Qe^{ \lambda_{i}}= M(1)\otimes \mathbb{C}[Q]e^{ \lambda_{i}}, \;
\; \; i=1,\dots, n$$ Given a Lie algebra element $a \otimes t^m \in
\widehat{\goth{sl}(n+1)}$, where $a \in \goth{sl}(n+1), m \in
\mathbb{Z}$, we will denote its action on an
$\widehat{\goth{sl}(n+1)}$-module using the notation $a(m)$. In particular, for $h \in
\goth{h}$ and $m \in \mathbb{Z}$, we have the operators $h(m)$ on
$V_P$:
$$h(0)(v\otimes \iota(e_\alpha)) = \langle h,\alpha \rangle (v\otimes \iota(e_\alpha)) $$ 
$$h(m) (v\otimes \iota(e_\alpha)) =( h(m) v\otimes \iota(e_\alpha)),$$
making $V_P$, $V_Q$, $V_Qe^{ \lambda_{1}}, \dots, V_Qe^{ \lambda_{n}}$ into 
$\widehat{\mathfrak{h}}$-modules.
This action may be extended using vertex operators to make $V_P$, $V_Q$, and $V_Qe^{\lambda_i}$ into
$\widehat{\goth{sl}(n+1)}$-modules (\cite{FLM}, \cite{FZ}, cf. \cite{LL}) with 
vertex operators satisfying
$$
Y(1 \otimes e^{\alpha},x) = \sum_{n \in \mathbb{Z}} x_{\alpha}(n)x^{-n-1}.
$$ 
In fact, we have that $L(\Lambda_i) \simeq V_Qe^{\lambda_i}$ for $i=0,\dots,n$.

We have natural operators 
\begin{eqnarray*}
e_{\alpha} \cdot e^{\beta}&=&
{\epsilon(\alpha, \beta)}e^{\alpha + \beta} 
\end{eqnarray*}
for $\alpha, \beta \in P,$  and that 
$$
x_\alpha(m) e_\lambda = c(\alpha,\lambda)e_\lambda x_\alpha(m + \langle \alpha, \lambda \rangle)
$$
for each $\alpha \in \Delta$ and $m \in \mathbb{Z}$, where 
$$
\epsilon: P \times P \rightarrow \mathbb{Z}/2(n+1)^2\mathbb{Z}
$$
is defined using 
a cocycle 
and 
$$
c: P \times P \rightarrow \mathbb{Z}/2(n+1)^2\mathbb{Z}
$$
is defined using a commutator map, both of which are 
 obtained by taking a central extension of $P$ by a cyclic group of order $2(n+1)^2$
(cf. \cite{FLM}, \cite{LL}).

As in \cite{S}, we need certain intertwining
operators among level $1$ standard modules.
We have intertwining operators
\begin{eqnarray} \label{IntwOp}{ \cal Y}( \cdot , x): L(\Lambda_r) &
 \longrightarrow & \mbox{Hom}(L(\Lambda_s), L( \Lambda_p))\{ x \} \\ w
 & \mapsto & {\cal Y}(w, x)=Y(w, x)e^{i\pi\lambda_r}c(\cdot,
 \lambda_r) \nonumber \end{eqnarray}
of type 
\begin{equation} 
\label{type} \left( \begin{array} {c} L(
 \Lambda_p) \\ \begin{array}{cc} L( \Lambda_r) & L( \Lambda_s)
 \end{array} \end{array} \right)  
\end{equation}
if and only if $p \equiv r+s \;
\mbox{mod}\; (n+1)$ (cf. \cite{DL}).
Given such an intertwining operator, we define a map
$${\cal Y}_c(e^{\lambda_r}, x): L(\Lambda_s) \longrightarrow
L(\Lambda_p)$$ by
$${\cal Y}_c(e^{\lambda_r}, x) = \mathrm{Res}_x x^{-1-\langle
  \lambda_r , \lambda_s\rangle}{\cal Y}(e^{\lambda_r}, x).$$ For each
$\alpha \in \Delta_+$, we have that
\begin{equation} \label{calYcommute}
[Y(e^{\alpha}, x_1), {\cal Y}_c(e^{\lambda_r}, x_2)]=0,
\end{equation}
which implies
\begin{equation} 
[x_\alpha(m), {\cal Y}_c(e^{\lambda_r}, x_2)]=0
\end{equation}
for each $m \in \mathbb{Z}$, a consequence of the Jacobi identity for
intertwining operators among standard modules.

We extend the operators $e_\lambda$, $\lambda \in P$, to operators 
on $V_p^{\otimes k}$, $k$ a positive integer, by defining:
$$
e_{\lambda}^{\otimes k} = e_\lambda \otimes \cdots \otimes e_\lambda:
V_P^{\otimes k} \rightarrow V_P^{\otimes k}.
$$

For any standard $\widehat{\goth{sl}(n+1)}$-module
$L(\Lambda)$ of positive integral level $k$, its highest weight $\Lambda $ is of the form
$$\Lambda = k_0\Lambda_0 + \dots + k_n \Lambda_n$$ for some
nonnegative integers $k_0,\dots,k_n$ satisfying $k_0 + \dots + k_n =
k$. Consider the space
 \begin{equation}
V_P^{\otimes k}= \underbrace{V_P \otimes \cdots \otimes V_P}_{k \; \;\mbox{times}},
\end{equation} 
and let 
\begin{equation}
v_{i_1, \dots, i_k}= v_{\Lambda_{i_1}} \otimes \cdots \otimes
v_{\Lambda_{i_k}} \in V_P^{\otimes k},
\end{equation}
where exactly $k_i$ indices are equal to $i$ for each
$i=0,\dots,n$. Then, we have that
$v_{i_1, \dots, i_k}$ is a highest weight vector for
$\widehat{\goth{sl}(n+1)}$, and
\begin{equation}
L(\Lambda) \simeq U(\widehat{{\goth{sl}(n+1)}}) \cdot v_{i_1, \dots,
i_k} \subset V_P ^{\otimes k}
\end{equation}
(cf. \cite{K}).
 We also
have a natural vertex operator algebra structure on $L(k\Lambda_0)$
and $L(k\Lambda_0)$-module structure on $L(\Lambda)$ given by:

\begin{theorem}(\cite{FZ}, \cite{DL}, \cite{Li1}; cf. \cite{LL}) 
The standard module $L(k \Lambda_0)$ has a natural vertex operator
algebra structure. The level $k$ standard
$\widehat{{\goth{sl}(n+1)}}$-modules provide a complete list of
irreducible $L(k \Lambda_0)$-modules. 
\end{theorem}

Let $\omega$ denote the Virasoro vector in $L(k \Lambda_0)$. We have a
natural representation of the Virasoro algebra on each $L(\Lambda)$
given by
\begin{equation}\label{Yomega}
Y_{L(\Lambda)}(\omega,x)=\sum_{m \in \mathbb{Z}} L(m)x^{-m-2}
\end{equation}
The operators $L(0)$ defined in (\ref{Yomega}) provide each
$L(\Lambda)$ of level $k$ with a grading, which we refer to as the
{\em weight} grading:
\begin{equation} \label{WeightGrading}
L(\Lambda) = \coprod_{s \in \mathbb{Z}} L(\Lambda)_{(s+h_{\Lambda})}
\end{equation}
where $h_\Lambda \in \mathbb{Q}$ and depends on $\Lambda$. In
particular, we have the grading
\begin{equation}
L(k\Lambda_0) = \coprod_{s \in \mathbb{Z}} L(\Lambda)_{(s)}.
\end{equation}
We denote the weight of an element $a \cdot v_{\Lambda} \in W(\Lambda)$ by 
$\mathrm{wt}(a\cdot v_{\Lambda})$. We will also write 
$$\mathrm{wt} (x_{\alpha}(m)) = -m,$$
where we view $x_\alpha(m)$ both as an operator and as an element of 
$U(\bar{\goth{n}})$.

We also have $n$ distinct {\em charge} gradings on each $L(\Lambda)$
of level $k$, given by the eigenvalues of the operators $\lambda_i(0)$
for $i=1,\dots,n$:
\begin{equation} \label{ChargeGrading}
L(\Lambda) = \coprod_{r_i \in \mathbb{Z}} L(\Lambda)_{[r_i+\langle\lambda_i,\Lambda\rangle]}.
\end{equation}
We call these the $\lambda_i$-{\em charge} gradings. An element of
$L(\Lambda)$ with $\lambda_i$-charges $n_i$ for $i=1,\dots,n$ 
has {\em total charge} $\sum_{i=1}^n n_i$. The gradings
(\ref{WeightGrading}) and (\ref{ChargeGrading}) are compatible, and we
have that
\begin{equation} \label{TripleGrading}
L(\Lambda) = \coprod_{r_1,\dots,r_n,s \in \mathbb{Z}}
L(\Lambda)_{r_1+\langle\lambda_1,\Lambda\rangle,\dots,
  r_n+\langle\lambda_n\Lambda\rangle;s+h_{\Lambda}}.
\end{equation}

\section{Principal subspaces and ideals}
\setcounter{equation}{0}
We are now ready to define our main object of study.
Consider the $\widehat{\goth{sl}(n+1)}$-subalgebra 
\begin{equation} 
\bar{\goth{n}}= \goth{n} \otimes \mathbb{C}[t, t^{-1}].
\end{equation}
The Lie algebra $\bar{\goth{n}}$ has the following
important subalgebras:
$$
\bar{\goth{n}}_{-} = \goth{n} \otimes t^{-1}\mathbb{C}[t^{-1}]
$$
and
$$
\bar{\goth{n}}_{+} = \goth{n} \otimes \mathbb{C}[t]
$$
Let $U(\bar{\goth{n}})$ be the universal 
enveloping algebra of $\bar{\goth{n}}$. 
We recall that $U(\bar{\goth{n}})$ has the decomposition:
\begin{equation} \label{decomp}
U(\bar{\goth{n}}) = U(\bar{\goth{n}}_-) \oplus U(\bar{\goth{n}})\bar{\goth{n}}_+.
\end{equation}
Given a $\widehat{\goth{sl}(n+1)}$-module $L(\Lambda)$ of positive
integral level $k$ with highest weight vector $v_\Lambda$, the {\it
  principal subspace} of $L(\Lambda)$ is defined by:
$$W(\Lambda) = U(\bar{\goth{n}})\cdot v_\Lambda.$$
$W(\Lambda)$ inherits the grading (\ref{TripleGrading}), and we have that
\begin{equation} \label{WTripleGrading}
W(\Lambda) = \coprod_{r_1,\dots,r_n,s \in \mathbb{Z}}
W(\Lambda)_{r_1+\langle\lambda_1,\Lambda\rangle,\dots,
  r_n+\langle\lambda_n\Lambda\rangle;s+h_{\Lambda}}
\end{equation}
For convenience, we will use the notation
$$ W(\Lambda)'_{r_1,\dots,
  r_n;s}=W(\Lambda)_{r_1+\langle\lambda_1,\Lambda\rangle,\dots,
  r_n+\langle\lambda_n\Lambda\rangle;s+h_{\Lambda}}
$$
As in \cite{CLM1}-\cite{CLM2}, \cite{CalLM3}, \cite{C1}-\cite{C2},
 define the {\em multigraded dimension} of $W(\Lambda)$ by:
$$ \chi_{W(\Lambda)}(x_1,\dots,x_n,
 q)=tr_{W(\Lambda)}x_1^{\lambda_1}\cdots x_n^{\lambda_n} q^{L(0)}
$$
and its modification
$$
\chi_{W(\Lambda)}'(x_1,\dots,x_n, q)= x^{-\langle \lambda_1, \Lambda \rangle} \dots
x^{-\langle \lambda_n, \Lambda \rangle} q^{-h_\Lambda}\chi_{W(\Lambda)}(x_1,\dots,x_n, q) \in 
\mathbb{C}[[x_1, \dots x_n, q]]
$$ In particular, we have that
$$ \chi_{W(\Lambda)}'(x_1,\dots,x_n, q) = \sum_{r_1,\dots,r_n,s \in
  \mathbb{N}}\mathrm{dim}(W(\Lambda)'_{r_1,\dots,r_n;s})x^{r_1}\cdots
x^{r_n}q^s.
$$

For each such $\Lambda$, we have a surjective map
\begin{eqnarray} \label{surj1}
F_{\Lambda}: U(\widehat{\goth{g}}) & \longrightarrow &
L(\Lambda) \\ a &\mapsto& a \cdot v_{\Lambda} \nonumber
\end{eqnarray}
and its surjective restriction $f_{\Lambda}$:
\begin{eqnarray} \label{surj2}
f_{\Lambda}: U(\bar{\goth{n}}) & \longrightarrow &
W(\Lambda)\\ a & \mapsto & a \cdot v_{\Lambda}. \nonumber
\end{eqnarray}
A precise description of the kernels $\mbox{Ker}f_{\Lambda}$ for every
each $\Lambda =\sum_{i=0}^n k_i\Lambda_i$ gives a presentation of the
principal subspaces $W(\Lambda)$ for $\widehat{\goth{sl}(n+1)}$, as we
will now discuss.

Define the formal sums
\begin{equation} 
R_{t}^i= \sum_{m_1+\dots
  +m_n=-t}x_{\alpha_i}(m_1)x_{\alpha_i}(m_2)\cdots
x_{\alpha_i}(m_{k+1})
\end{equation}
and their truncations
\begin{equation} \label{r}
R_{M, t}^i= \sum_{\begin{array} {c}m_1+ \cdots + m_{k+1}=-t ,
    \\  m_1, \dots, m_{k+1} \leq M \end{array}}x_{\alpha_i}(m_1)\cdots
x_{\alpha_i}(m_{k+1})
\end{equation}
for $t \in \mathbb{Z}$, $M \in \mathbb{Z}$ and
$i=1,\dots,n$. 
Note that each $R_{M, t}^i \in U(\bar{\goth{n}})$ and 
the infinite sum $R_{t}^i \notin U(\bar{\goth{n}})$,
but $R_t^i$ is still well-defined as an operator on $W(\Lambda)$, since,
when acting on any element of $W(\Lambda)$, only finitely many of its terms
are nonzero.
Let $J$ be the left ideal of $U(\bar{\goth{n}})$
generated by the elements $R_{-1, t}^i$ for $t \geq k+1$ and $i=1,2$:
\begin{equation}
J = \sum_{i=1}^n \sum_{t \geq k+1} U(\bar{\goth{n}}) R_{-1,t}^i.
\end{equation}
Define a left ideal of $U(\bar{\goth{n}})$ by:
\begin{equation*} 
I_{k\Lambda_{0}}= J+U(\bar{\goth{n}})\bar{\goth{n}}_{+}
\end{equation*}
and for each $\Lambda = \sum_{i=0}^n k_i \Lambda_i$, define
$$I_\Lambda = I_{k\Lambda_0} + \sum_{\alpha \in \Delta_+}
U(\bar{\goth{n}})x_{\alpha}(-1)^{k+1-\langle \alpha, \Lambda \rangle
}.$$

\begin{conjecture}\label{presentations}
For each $\Lambda = k_0\Lambda_0 + \dots + k_n\Lambda_n$ with 
$k_0,\dots ,k_n, k \in \mathbb{N}$, $k \ge 1$, and $k_0 + \dots +k_n = k$, we have that
$$\mathrm{Ker}f_\Lambda = I_\Lambda.$$
For $\Lambda = k_0 \Lambda_0 + k_i \Lambda_i$,we have
\begin{equation} \label{suggested}
\mathrm{Ker}f_{k_0\Lambda_0 + k_i\Lambda_i} = I_{k \Lambda_0} + U(\bar{\goth{n}})x_{\alpha_i}(-1)^{k_0+1}.
\end{equation}
\end{conjecture}

In the cases that $\mathfrak{g} = \mathfrak{sl}(2)$ and $k\ge1$ or $\mathfrak{g}$
is of type $ADE$ and $k=1$, Conjecture \ref{presentations} has been proved in \cite{CalLM1}-\cite{CalLM3}. 
In the case that $\mathfrak{g} = \mathfrak{sl}(3)$
and $k\ge 1$ this conjecture has been proved in \cite{S}.
The presentations
(\ref{suggested}) are suggested by the bases found in \cite{G}, but an a
priori proof is lacking. The proof of (\ref{suggested}) will appear in future work.

We now give a partial proof of this conjecture, using the quasiparticle bases 
for principal subspaces obtained in \cite{G} to obtain the second term in (\ref{suggested})
for the cases $i=1$ and $i=n$.
As in \cite{G}, for each $i=1,\dots,n$, we define
${n}_{\alpha_i} = \mathbb{C} x_{\alpha_i}$ and 
$\bar{n}_{\alpha_i} = {n}_{\alpha_i} \otimes \mathbb{C}[t,t^{-1}]$.
For each $i=1,\dots,n$, define the operators
\begin{eqnarray*}
 x_{M\alpha_i}(m) = \sum_{m_1 + \dots + m_M = m} x_{\alpha_i}(m_1)\cdots x_{\alpha_i}(m_M),
\end{eqnarray*}
called  quasiparticles of color $i$ and charge $M$ in \cite{G}.
These operators act in a well defined way on any element of $W(\Lambda)$, $\Lambda$ a dominant 
integral weight, in the sense that, when applied to an element of $W(\Lambda)$, only finitely
many terms are nonzero. We will consider $\Lambda$ of the form
$$
\Lambda = k_0 \Lambda_0 + k_j \Lambda_j
$$
for $j=1,\dots,n$ and $k_0,k_j \in \mathbb{N}$ with $k_0 + k_j = k$.
Define
$$
j_{t} :=  \left\{ \begin{array}{l} \ 0\;\; \mbox{for} \;
0 < t \leq  k_{0} \\  j\;\;  \mbox{for} \;
k_{0}< t \leq k = k_{0}+ k_{j}.\end{array}
\right.\label{4.1}
$$
In \cite{G}, Georgiev proved that
$$
W(\Lambda) = U(\bar{\mathfrak{n}}_{\alpha_{n}}) \cdots
U(\bar{\mathfrak{n}}_{\alpha_{1}}) \cdot v_\Lambda
$$
and also that the set of operators
\begin{eqnarray}\label{bases}
 \lefteqn{\mathfrak{B}_{W(\Lambda)} :=}
\end{eqnarray}
$$
\bigsqcup_{\begin{array}{c}{\scriptstyle 0 \leq
n_{r_{n}^{(1)},n} \leq \cdots \leq  n_{1,n}\leq k}\\{\scriptstyle
\cdots\cdots\cdots}\\{\scriptstyle 0 \leq
n_{r_{1}^{(1)},1} \leq \cdots \leq n_{1,1}\leq
k}\end{array}}
\left\{x_{n_{r_{n}^{(1)},n}\alpha_{n}}(m_{r_{n}^{(1)},n})
\cdots x_{n_{1,n}\alpha_{n}}(m_{1,n}) \cdots \cdots
x_{n_{r_{1}^{(1)},1}\alpha_{1}}(m_{r_{1}^{(1)},1}) \cdots
x_{n_{1,1}\alpha_{1}}(m_{1,1}) \begin{array}{c}\\ \\ \\
\end{array}\right|$$
$$\left|\begin{array}{ll} m_{p,i} \in \mathbb{Z},\; 1\leq i \leq n, \; 1 \leq
p \leq r_{i}^{(1)};&\\m_{p,i} \leq
\sum_{q=1}^{r_{i-1}^{(1)}} \mbox{min}\,\{n_{p,i}, n_{q,i-1}\}  -
\sum_{t=1}^{n_{p,i}}\delta_{i,j_{t}}
-\sum_{p>p'>0} 2\mbox{min}\, \{n_{p,i}, n_{p',i}\} -
n_{p,i};&\\m_{p+1,i} \leq m_{p,i} -2n_{p,i}\;\;\mbox{for}\;\;
n_{p+1,i} = n_{p,i}& \end{array}\right\},
$$
where $r_{0}^{(1)} := 0$, forms a basis for $W(\Lambda)$ when 
applied to $v_\Lambda$. 
It is important to notice that, for each $j = 1,\dots ,n$, 
$$
\mathfrak{B}_{W(k\Lambda_j)} \subset \mathfrak{B}_{W(\Lambda_0 + (k-1)\Lambda_j)}
\subset \cdots \subset \mathfrak{B}_{W(k\Lambda_0)}
$$

Consider the maps
$$
f_{k\Lambda_0}: U(\bar{\mathfrak{n}}) \longrightarrow W(k\Lambda_0)
$$
and
$$
1^{\otimes k_0} \otimes \mathcal{Y}_c(e^{\lambda_j},x)^{k_j}: 
W(k\Lambda_0) \longrightarrow W(k_0\Lambda_0 + k_j\Lambda_j).
$$
Composing these maps, we obtain (after multiplication by a scalar):
$$
f_{k_0\Lambda_0 + k_j \Lambda_j} = (1^{\otimes k_0} \otimes \mathcal{Y}_c(e^{\lambda_j},x)^{k_j}) \circ 
f_{k\Lambda_0},
$$
which implies that 
\begin{eqnarray}\label{kercontain}
\mbox{Ker}f_{k\Lambda_0} \subset \mbox{Ker}f_{k_0\Lambda_0 + k_j \Lambda_j} 
\end{eqnarray}
Using these maps, along with the bases above, we obtain the following proposition:
\begin{proposition} Fix $k_0, k_1, k_n, k \in \mathbb{N}$ with $k_0 + k_1 = k_0 + k_n = k >0$. Then
\begin{eqnarray}\label{1discrepancy}
\mathrm{Ker}f_{k_0\Lambda_0 + k_1 \Lambda_1} = 
\mathrm{Ker}f_{k\Lambda_0}+ U(\bar{\goth{n}})x_{\alpha_1}(-1)^{k_0 + 1}
\end{eqnarray}
and
\begin{eqnarray}\label{ndiscrepancy}
\mathrm{Ker}f_{k_0\Lambda_0 + k_n \Lambda_n} = 
 \mathrm{Ker}f_{k\Lambda_0} + U(\bar{\goth{n}})x_{\alpha_n}(-1)^{k_0 + 1}.
\end{eqnarray}
\end{proposition}

{\em Proof:} We prove (\ref{1discrepancy}). Note that (\ref{ndiscrepancy})
will follow by an identical proof when the bases (\ref{bases}) are rewritten 
as subsets of $U(\bar{\mathfrak{n}}_{\alpha_{1}}) \cdots
U(\bar{\mathfrak{n}}_{\alpha_{n}})$.

First, the fact that 
$$
U(\bar{\goth{n}})x_{\alpha_1}(-1)^{k_0 + 1} \subset \mbox{Ker}f_{k_0\Lambda_0 + k_1 \Lambda_1}
$$
follows immediately from the fact that 
$$
x_{\alpha_1}(-1)^2 \cdot v_{\Lambda_0} = 0
$$
and
$$
x_{\alpha_1}(-1)\cdot v_{\Lambda_1} = 0,
$$
and, along with (\ref{kercontain}) gives us that 
$$
\mbox{Ker}f_{k\Lambda_0}+ U(\bar{\goth{n}})x_{\alpha_1}(-1)^{k_0 + 1} \subset
\mbox{Ker}f_{k_0\Lambda_0 + k_1 \Lambda_1} 
$$

We now show the reverse inclusion. Using the spanning argument for (\ref{bases}) in \cite{G}, we may
write 
$a\cdot v_{k\Lambda_0} \in \mathrm{span}(\mathfrak{B}_{W(k\Lambda_0)}\cdot v_{k\Lambda_0})$ as
\begin{eqnarray}\label{rewriting}
a\cdot v_{k\Lambda_0} = (b + cx_{\alpha_1}(-1)^{k_0 + 1})\cdot v_{k\Lambda_0}
\end{eqnarray}
for some $b \in \mathrm{span}({\mathfrak{B}_{W(k_0\Lambda_0 + k_1\Lambda_1)}})
\subset \mathrm{span}({\mathfrak{B}_{W(k\Lambda_0)}})$ and $c \in U(\bar{\mathfrak{n}}_{\alpha_{n}}) \cdots
U(\bar{\mathfrak{n}}_{\alpha_{1}})$. Suppose that  $a \in \mbox{Ker}f_{k_0\Lambda_0 + k_1 \Lambda_1} $.
If $a \in \mbox{Ker}f_{k\Lambda_0}$, we are done. So, suppose that $a \cdot v_{k\Lambda_0} \neq 0$.
Then, applying the map $(1^{\otimes k_0} \otimes \mathcal{Y}_c(e^{\lambda_1},x)^{k_1})$
and using (\ref{rewriting}) we have 
\begin{eqnarray*}
0 &=& a\cdot v_{k_0\Lambda_0 + k_1 \Lambda_1}\\
&=& (b + cx_{\alpha_1}(-1)^{k_0 + 1})
\cdot v_{k_0\Lambda_0 + k_1 \Lambda_1}\\
&=& b \cdot v_{k_0\Lambda_0 + k_1 \Lambda_1},
\end{eqnarray*}
which implies that $b=0$, so that $a \cdot v_{k\Lambda_0} = cx_{\alpha_1}(-1)^{k_0 + 1}\cdot v_{k\Lambda_0}$,
and $a \in \mbox{Ker}f_{k\Lambda_0} +  U(\bar{\mathfrak{n}})x_{\alpha_1}(-1)^{k_0 + 1}$.

\begin{remark}\em
 As an alternative to the spanning argument in \cite{G}, we may obtain
 (\ref{rewriting}) by solving a linear system of equations. 
 To obtain this system of equations, it suffices to rewrite only elements 
 of $\mathfrak{B}_{W(k\Lambda_0)} \setminus \mathfrak{B}_{W(k_0\Lambda_0 + k_1\Lambda_1)}$.
Given such an element $a \in \mathfrak{B}_{W(k\Lambda_0)} \setminus \mathfrak{B}_{W(k_0\Lambda_0 + k_1\Lambda_1)}$,
expand all the terms and
look at all summands in $a \cdot v_{k\Lambda_0}$
 not in $U(\bar{\mathfrak{n}}_{\alpha_{n}}) \cdots
U(\bar{\mathfrak{n}}_{\alpha_{1}})x_{\alpha_1}(-1)^{k_0+1}\cdot v_{k\Lambda_0}$
and match them up with linear combinations of expanded elements from 
  $\mathfrak{B}_{W(k_0\Lambda_0 + k_1\Lambda_1)}\cdot v_{k\Lambda_0}$ with
 the same charges and weight.
The term $cx_{\alpha_1}(-1)^{k_0+1}$
 arises from adding in terms from $U(\bar{\mathfrak{n}}_{\alpha_{n}}) \cdots
U(\bar{\mathfrak{n}}_{\alpha_{1}})x_{\alpha_1}(-1)^{k_0+1}\cdot v_{k\Lambda_0}$ to make both sides
of (\ref{rewriting}) equal.
 Take, for example, the operator
 $$
 x_{2\alpha_1}(-6)x_{3\alpha_1}(-5) \in \mathfrak{B}_{W(k\Lambda_0)},
 $$
 where we assume $k>5$. We have, by definition, that
 $$
  x_{2\alpha_1}(-6)x_{3\alpha_1}(-5) \notin \mathfrak{B}_{W(\Lambda_0 + (k-1)\Lambda_1)}.
 $$
 Applying this operator to $v_{\Lambda_0 + (k-1)\Lambda_1}$, we have that
 $$x_{2\alpha_1}(-6)x_{3\alpha}(-5) \cdot v_{\Lambda_0 + (k-1)\Lambda_1} \neq 0.$$
 We may write $x_{2\alpha_1}(-6)x_{3\alpha_1}(-5)\cdot v_{k\Lambda_0}$ as a linear combination of
 $$
 x_{5\alpha_1}(-11)\cdot v_{k\Lambda_0},
 \  x_{\alpha_1}(-3)x_{4\alpha_1}(-8)\cdot v_{k\Lambda_0},
 \  x_{\alpha_1}(-4)x_{4\alpha_1}(-7)\cdot v_{k\Lambda_0}
 \in \mathfrak{B}_{W(\Lambda_0 + (k-1)\Lambda_1)}\cdot v_{k\Lambda_0}
 $$
 and an element of the form $U(\bar{\mathfrak{n}}_{\alpha_{n}}) \cdots
U(\bar{\mathfrak{n}}_{\alpha_{1}})x_{\alpha_1}(-1)^2 \cdot v_{k\Lambda_0}$, giving 
 precisely an expression of the form (\ref{rewriting}), as follows:
 \begin{eqnarray*}
  \lefteqn{x_{2\alpha_1}(-6)x_{3\alpha_1}(-5) \cdot v_{k\Lambda_0}}\\
  &&= (\frac{1}{2}x_{\alpha_1}(-3)x_{4\alpha_1}(-8) - 
  \frac{1}{10} x_{5\alpha_1}(-11) + 2x_{\alpha_1}(-4)x_{4\alpha_1}(-7)) \cdot v_{k\Lambda_0} 
  + cx_{\alpha_1}(-1)^2\cdot v_{k\Lambda_0} 
 \end{eqnarray*}
for some $c \in U(\bar{\mathfrak{n}}_{\alpha_{n}}) \cdots
U(\bar{\mathfrak{n}}_{\alpha_{1}})$.

 \em
\end{remark}

\section{A reformulation of the presentation problem}
\setcounter{equation}{0}
In this section we reformulate Conjecture \ref{presentations}, along
with all known presentations of principal subspaces, in terms of a natural completion of
$U(\bar{\goth{n}})$, which we denote by
$\widetilde{U(\bar{\goth{n}})}$. A version of this completion was
constructed in \cite{LW}, and we recall this construction, suitably adapted
to our present setting, in the appendix. In this section
only, for $\alpha \in \Delta$ and 
$n \in \mathbb{Z}$, will use the notation 
$x_\alpha(n)$ for completion elements $X_\alpha(n)$ from the appendix,
and no confusion should arise.

We may define a natural ``lifting" of the maps $f_\Lambda$:
\begin{eqnarray} 
\widetilde{f_{\Lambda}}: \widetilde{U(\bar{\goth{n}})} & \longrightarrow &
W(\Lambda)\\ a & \mapsto & a \cdot v_{\Lambda}. \nonumber
\end{eqnarray}
Indeed, given $a \in \widetilde{U(\bar{\goth{n}})}$, we may uniquely
express $a$ as $a = b+c$ for some $b \in U(\bar{\goth{n}}_-)$ and $c
\in \widetilde{U(\bar{\goth{n}})\bar{\goth{n}}_+}$ (by (\ref{compdecomp})), and define $a
\cdot v_\Lambda = b \cdot v_\Lambda$. That is, we let $c$ act as $0$.

We now reformulate Conjecture \ref{presentations} in terms of finding
$\mbox{Ker}\widetilde{f_{\Lambda}}$.  Recall the formal sums
$$
R_t^i = \sum_{m_1 + \dots + m_{k+1} = -t} x_{\alpha_i}(m_1)\cdots x_{\alpha_i}(m_{k+1}),
$$
which are well defined as operators on each $W(\Lambda)$.
It is important to note that each $R_t^i$ is not an element of
$\widetilde{U(\bar{\goth{n}})}$, so we seek natural representatives
for $R_t^i$ in $\widetilde{U(\bar{\goth{n}})}$, in the sense that, when viewed as operators
on $W(\Lambda)$, these representatives are equal to $R_t^i$.

Let $\mathcal{A}$ denote the set of finite sequences of integers.
Given a sequence of integers $A = (m_1,\dots,m_{k+1}) \in
\mathcal{A}$, define a function
\begin{eqnarray} 
\#: \mathbb{Z} \times A & \longrightarrow & \mathbb{N}\\ (n,A) &
\mapsto & \mathrm{number\ of\ occurrences\ of\ }\ n\ \mathrm{ \ in}
\ A. \nonumber
\end{eqnarray}
For any sequence in $(m_1,...,m_{k+1}) \in \mathcal{A}$, define
$$
A_{m_1,\dots,m_{k+1}} = \{\#(n,A)| n \in \mathbb{Z} \}\setminus \{0\} = \{n_1,\dots,n_j\}
$$
where $n_1,\dots,n_j$ are positive integers and $n_1 + \dots + n_j = k+1$.
Define integers
$$
c_{m_1,...,m_{k+1}} = \binom{k+1}{n_1,...,n_j} = \frac{(k+1)!}{(n_1)!\dots (n_j)!}.
$$
We define
\begin{equation}\label{-1rep}
\mathcal{R}_t^i = R_{-1,t}^i + \sum_{\begin{array} {c} m_1 \le \dots
    \le m_{k+1}, \\ m_1+ \cdots + m_{k+1}=-t, \\ m_{k+1} \ge
    0 \end{array}}c_{m_1, \dots, m_{k+1}} x_{\alpha_i}(m_1) \cdots
x_{\alpha_i}(m_{k+1}),
\end{equation}
which is clearly in $\widetilde{U(\bar{\goth{n}})}$. We may also
write, for each $\mathcal{R}_t^i$,
\begin{equation}\label{Mrep}
\mathcal{R}_t^i = R_{M,t}^i + \sum_{\begin{array} {c} m_1 \le \dots
    \le m_{k+1}, \\ m_1+ \cdots + m_{k+1}=-t, \\ m_{k+1} \ge
    M+1 \end{array}}c_{m_1, \dots, m_{k+1}} x_{\alpha_i}(m_1) \cdots
x_{\alpha_i}(m_{k+1}),
\end{equation}
and, as elements of $\widetilde{U(\bar{\goth{n}})}$, (\ref{-1rep}) and
(\ref{Mrep}) are equal.
\begin{remark}\em
As mentioned above, the formal sums
$$
R_t^i = \sum_{m_1 + \dots + m_{k+1} = -t} x_{\alpha_i}(m_1)\cdots x_{\alpha_i}(m_{k+1}).
$$ are not elements of $\widetilde{U(\bar{\goth{n}})}$. Informally,
$R_t^i$ is in a sense a ``limit" of (\ref{Mrep}), i.e.
$$
R_t^i = \lim_{M \rightarrow \infty}\bigg( R_{M,t}^i + \sum_{\begin{array} {c} m_1 \le \dots \le m_{k+1}, \\ m_1+ \cdots + m_{k+1}=-t,
\\ m_{k+1} \ge M+1 \end{array}}c_{m_1, \dots, m_{k+1}} x_{\alpha_i}(m_1)
\cdots x_{\alpha_i}(m_{k+1})\bigg),
$$ where infinitely many relations in $\widetilde{I}$ need to be
applied to obtain $R_t^i$ from $\mathcal{R}_t^i$. 
However, as operators on $W(\Lambda)$, $\mathcal{R}_t^i$ and $R_t^i$ are equal.\em
\end{remark}

\begin{lemma}
Let $\alpha \in \Delta_+$ and $m \in \mathbb{N}$. Then, for any
$i=1,\dots,n$ and $t \in \mathbb{Z}$ we have that
$$
\mathcal{R}_t^i x_{\alpha}(-m) = x_\alpha(-m)\mathcal{R}_t^i + x_\alpha(0)R_{t+m}^i + c
$$ for some $c \in \widetilde{U(\bar{\goth{n}})\bar{\goth{n}}_+}$. In
particular,
$$
\mathcal{R}_t^i x_{\alpha}(-m) \in I_{k\Lambda_0} + \widetilde{U(\bar{\goth{n}})}\bar{\goth{n}}_+.
$$
\end{lemma}

{\em Proof:} First, suppose that $\alpha + \alpha_i \in \Delta_+$. We may write
$$
\mathcal{R}_t^i = R_{m,t}^i + \sum_{\begin{array} {c} m_1 \le \dots \le m_{k+1}, \\ m_1+ \cdots + m_{k+1}=-t,
\\ m_{k+1} \ge m+1 \end{array}}c_{m_1, \dots, m_{k+1}} x_{\alpha_i}(m_1)
\cdots x_{\alpha_i}(m_{k+1}).
$$
By definition of $\widetilde{U(\bar{\goth{n}})\bar{\goth{n}}_+}$,
$$
\sum_{\begin{array} {c} m_1 \le \dots \le m_{k+1}, \\ m_1+ \cdots + m_{k+1}=-t,
\\ m_{k+1} \ge m+1 \end{array}}c_{m_1, \dots, m_{k+1}}x_{\alpha_i}(m_1)
\cdots x_{\alpha_i}(m_{k+1}) x_\alpha(-m) \in 
 \widetilde{U(\bar{\goth{n}})\bar{\goth{n}}_+}.
$$
For $R_{m,t}^ix_{\alpha}(-m)$, we may write
\begin{eqnarray*}
R_{m,t}^ix_{\alpha}(-m) &=& \sum_{\begin{array} {c} m_1+ \cdots +
    m_{k+1}=-t, \\ m_1,\dots , m_{k+1} \le m \end{array}}
x_{\alpha_i}(m_1)\cdots x_{\alpha_i} (m_{k+1}) x_\alpha(-m)\\ &=&
\sum_{j=1}^{k+1}\sum_{\begin{array} {c} m_1+ \cdots + m_{k+1}=-t,
    \\ m_1,\dots , m_{k+1} \le m \end{array}} C_{\alpha_i,\alpha}
x_{\alpha_i}(m_1)\cdots x_{\alpha_i + \alpha}(m_j-m)\cdots
x_{\alpha_i} (m_{k+1})\\ && + x_{\alpha}(-m)R_{m,t}^i\\ &=&
\sum_{j=1}^{k+1}\sum_{\begin{array} {c} m_1+ \cdots + m_{k+1}=-t-m,
    \\ m_1,\dots , m_{k+1} \le m \end{array}} C_{\alpha_i,\alpha}
x_{\alpha_i}(m_1)\cdots x_{\alpha_i + \alpha}(m_j)\cdots x_{\alpha_i}
(m_{k+1})\\ && + b + x_{\alpha}(-m)R_{m,t}^i
\end{eqnarray*}
for some $b \in U(\bar{\goth{n}})\bar{\goth{n}}_+$. We have that
\begin{eqnarray*}
\sum_{j=1}^{k+1}\sum_{\begin{array} {c} m_1+ \cdots + m_{k+1}=-t-m,
    \\ m_1,\dots , m_{k+1} \le m \end{array}} C_{\alpha_i,\alpha}
x_{\alpha_i}(m_1)\cdots x_{\alpha_i + \alpha}(m_j)\cdots x_{\alpha_i}
(m_{k+1}) + b + x_{\alpha}(-m)R_{m,t}^i\\ = [x_{\alpha}(0),R_{m+t}^i]
+ b + x_{\alpha}(-m)R_{m,t}^i,
\end{eqnarray*}
establishing our claim when $\alpha + \alpha_i \in \Delta_+$.  If
$\alpha + \alpha_i \notin \Delta_+$ the claim is clear
since
$$
\mathcal{R}_t^i x_\alpha(-m) = x_\alpha(-m)\mathcal{R}_t^i \in I_{k\Lambda_0} +
\widetilde{U(\bar{\goth{n}})\bar{\goth{n}}_+},
$$
concluding our proof.\\

Using an almost identical argument, we have that
\begin{corollary}\label{commcorr}
If $a \in U(\bar{\goth{n}}_-)$ and $t \in \mathbb{Z}$, we have that
$$
\mathcal{R}_t^i a \in I_{k\Lambda_0} + \widetilde{U(\bar{\goth{n}})\bar{\goth{n}}_+}.
$$
\end{corollary}
{\em Proof:} It suffices to show that the claim holds for monomials
$$
x_{\beta_1}(-m_1) \dots x_{\beta_j}(-m_j) \in
U(\bar{\goth{n}}).
$$
This follows immediately using the same argument as above, and writing
$$
\mathcal{R}_t^i = R_{m_1+\dots+m_j,t}^i + \sum_{\begin{array} {c} m_1 \le \dots \le m_{k+1}, \\ m_1+ \cdots + m_{k+1}=-t,
\\ m_{k+1} \ge m_1+\dots+m_j+1 \end{array}}c_{m_1, \dots, m_{k+1}} x_{\alpha_i}(m_1)
\cdots x_{\alpha_i}(m_{k+1}).
$$

As in \cite{C1}-\cite{C2} and \cite{CalLM3}, let $\mathcal{J}$ be the
two sided ideal of $\widetilde{U(\bar{\goth{n}})}$ generated by the
$\mathcal{R}_t^i$, $i=1,\dots,n$ and $t \ge k+1$.  As in
\cite{CalLM3}, we have the following theorem:
\begin{theorem}\label{compidealequivalence}
We may describe $I_{k\Lambda_0}$ by:
\begin{equation}
I_{k\Lambda_0} \equiv {\mathcal J} \; \; \;
\mbox{\rm modulo} \; \; \;
\widetilde{U(\bar{\goth{n}})\bar{\goth{n}}_{+}}.
\end{equation}
and moreover,
for  $I_\Lambda$, we have:
\begin{equation} \label{equation_2}
I_{\Lambda}\equiv {\mathcal J} + \sum_{\alpha \in \Delta_+}
U(\bar{\goth{n}})x_{\alpha}(-1)^{k+1-\langle \alpha , \Lambda \rangle}
\; \; \; \mbox{\rm modulo} \; \; \;
\widetilde{U(\bar{\goth{n}})\bar{\goth{n}}_{+}}.
\end{equation}
\end{theorem}
{\em Proof:} We first show that $$I_{k\Lambda_0} \subset {\mathcal J}
\ \mbox{\rm
  modulo}\ \widetilde{U(\bar{\goth{n}})\bar{\goth{n}}_{+}}.$$ Indeed,
any element $a \in I_{k\Lambda_0}$ may be written as
$$
a = \sum_{i=1}^n a_i R_{-1,t}^{i} + b
$$ for some $a_i \in U(\bar{\goth{n}})$ and $b \in
U(\bar{\goth{n}})\bar{\goth{n}}_+$. It suffices to show that each $a_i
R_{-1,t}^{i} \in \mathcal{J} +
\widetilde{U(\bar{\goth{n}})\bar{\goth{n}}_{+}}$. Indeed, 
we may write $R_{-1,t}^{i} = \mathcal{R}_t^i + c$ for some $c \in
\widetilde{U(\bar{\goth{n}})\bar{\goth{n}}_{+}}$, and we clearly have
that $a_iR_{-1,t}^{i} = a_iR_t^i + a_i c \in \mathcal{J} +
\widetilde{U(\bar{\goth{n}})\bar{\goth{n}}_{+}}$.

It remains to show that
$$
{\mathcal J} \subset I_{k\Lambda_0}
 \  \mbox{\rm modulo}\  \widetilde{U(\bar{\goth{n}})\bar{\goth{n}}_{+}}.
$$
It suffices to prove
$$
a \mathcal{R}_t^i b \in I_{k\Lambda_0} + \widetilde{U(\bar{\goth{n}})\bar{\goth{n}}_+}
$$ for all $a,b \in \widetilde{U(\bar{\goth{n}})}.$ By
(\ref{compdecomp}), we may write $b = b_1 + b_2$ for some $b_1 \in
U(\bar{\goth{n}}_-)$ and $b_2 \in
\widetilde{U(\bar{\goth{n}})\bar{\goth{n}}_+}$.  Clearly
$a\mathcal{R}_t^ib_2 \in
\widetilde{U(\bar{\goth{n}})\bar{\goth{n}}_+}$, and so
$$
a\mathcal{R}_t^i b \equiv a\mathcal{R}_t^ib_1 \ \mathrm{modulo} \ \widetilde{U(\bar{\goth{n}})\bar{\goth{n}}_+}.
$$
By Corollary \ref{commcorr}, we have that 
$$
\mathcal{R}_t^ib_1 \in I_{k\Lambda_0} + \widetilde{U(\bar{\goth{n}})\bar{\goth{n}}_+},
$$
so it suffices to show that 
$$
a R_{-1,t}^i \in I_{k\Lambda_0} + \widetilde{U(\bar{\goth{n}})\bar{\goth{n}}_+}.
$$

Using the notation from the appendix,
we have that $a = [\mu]$ for some $\mu \in F(\Delta_+)$, and we may write
$$ \mu = \sum_{c \in \mathrm{Supp}(\mu)}\mu(c)X(c) = \sum_{c \in
  \mathrm{Supp}_{t}(\mu)}\mu(c)X(c) + \sum_{c \in
  \mathrm{Supp}(\mu)\setminus\mathrm{Supp}_{t}(\mu)} \mu(c)X(c).
$$
The sum $\sum_{c \in \mathrm{Supp}_{t}(\mu)}\mu(c)X(c)$ is finite, so we have that 
$$
\sum_{c \in \mathrm{Supp}_{t}(\mu)}[\mu(c)X(c)]R_{-1,t}^i \in I_{k\Lambda_0},
$$ and by definition of
$\widetilde{U(\bar{\goth{n}})\bar{\goth{n}}_+}$ we have that
$$
\sum_{c \in \mathrm{Supp}(\mu)\setminus\mathrm{Supp}_{t}(\mu)} 
[\mu(c)X(c)]R_{-1,t}^i \in \widetilde{U(\bar{\goth{n}})\bar{\goth{n}}_+},
$$
establishing
$$I_{k\Lambda_0} \equiv {\mathcal J} \  \mbox{\rm modulo}\  \widetilde{U(\bar{\goth{n}})\bar{\goth{n}}_{+}}.$$
The fact that
$$
I_{\Lambda}\equiv {\mathcal J} +
\sum_{\alpha \in \Delta_+} U(\bar{\goth{n}})x_{\alpha}(-1)^{k+1-\langle \alpha , \Lambda \rangle} \; \; \; \mbox{\rm modulo} \; \; \;
\widetilde{U(\bar{\goth{n}})\bar{\goth{n}}_{+}}
$$ follows immediately, establishing our theorem.  \\

As a consequence of Theorem \ref{compidealequivalence}, along with the
results of \cite{CalLM1} - \cite{CalLM3} and \cite{S}, we have that:
\begin{theorem}
In the case where $\mathfrak{g} = \mathfrak{sl}(n+1)$ with:
\begin{itemize}
\item $n=1$ and $\Lambda = k_0\Lambda_0 + k_1\Lambda_1$ with $k_0 + k_1 = k\ge 1$
\item$n=2$ and $\Lambda = k_0 \Lambda_0 + k_1\Lambda_1 + k_2\Lambda_2$  with $k_0 + k_1 + k_2 = k\ge 1$
\item $n\ge 3$ and $\Lambda = \Lambda_i$ with $i=0,\dots,n$
\end{itemize}
or $\mathfrak{g}$ is of type $D$ or $E$  
with $k=1$
 we have that 
$$\mathrm{Ker}f_\Lambda \equiv \widetilde{I_\Lambda} \; \; \; \mbox{\rm
   modulo} \; \; \; \widetilde{U(\bar{\goth{n}})\bar{\goth{n}}_{+}}.$$
\end{theorem}

We reformulate Conjecture \ref{presentations} as follows:
\begin{conjecture}\label{compconjecture}
Suppose $\mathfrak{g} = \mathfrak{sl}(n+1)$,
$k_0,\dots,k_n,k \in \mathbb{N}$ with $k \ge 1$ and $k_0 + \dots + k_n
= k$.  For each $\Lambda = k_0\Lambda_0 + \dots + k_n \Lambda_n$, we
have that
$$\mathrm{Ker}f_\Lambda \equiv \widetilde{I_\Lambda}  \; \; \; \mbox{\rm modulo} \; \; \;
\widetilde{U(\bar{\goth{n}})\bar{\goth{n}}_{+}}$$
\end{conjecture}
or that 
\begin{conjecture}
In the context of Conjecture \ref{compconjecture},
for each $\Lambda = k_0\Lambda_0 + \dots + k_n \Lambda_n$, we have that
$$
\mathrm{Ker}\widetilde{f_\Lambda} = \widetilde{I_\Lambda}.
$$
\end{conjecture}

\section{Exact sequences and multigraded dimensions}
\setcounter{equation}{0}
In this section, we construct exact sequences among the principal
subspaces of certain standard modules, and use these to find
multigraded dimensions.

 Given $\lambda \in P$ and character $\nu : Q \longrightarrow
 \mathbb{C}^*$, we define a map $\tau_{\lambda,\nu}$ on
 $\bar{\goth{n}}$ by
\begin{equation*}
\tau_{\lambda, \nu}(x_{\alpha}(m))=\nu(\alpha) x_{\alpha}(m-\langle
\lambda, \alpha \rangle)
\end{equation*}
for $\alpha \in \Delta_{+}$ and $m \in \mathbb{Z}$.  It is easy to see
that $\tau_{\lambda,\nu}$ is an automorphism of $\bar{\goth{n}}$.  The
map $\tau_{\lambda,\nu}$ extends canonically to an automorphism of
$U(\bar{\goth{n}})$, also denoted by $\tau_{\lambda,\nu}$, given by
\begin{equation}\label{def-lambda_i} 
\tau_{\lambda,\nu}(x_{\beta_1}(m_1) \cdots x_{\beta_r}(m_r))=
\nu(\beta_1+\cdots + \beta_r)
x_{\beta_1}(m_1-\langle \lambda, \beta_1 \rangle ) \cdots
x_{\beta_r}(m_r-\langle \lambda, \beta_r \rangle)
\end{equation} 
for $\beta_1, \dots, \beta_r \in \Delta_{+}$ and $m_1, \dots, m_r \in
\mathbb{Z}$.
In particular, we have that
\begin{eqnarray}
e_\lambda^{\otimes k} (a\cdot v_\Lambda) =
\tau_{\lambda, c_{-\lambda}}(a)\cdot e_\lambda^{\otimes k}v_\Lambda
\end{eqnarray}
where $\lambda \in P$, $\Lambda$ is a dominant integral weight of
$\widehat{\mathfrak{sl}(n+1)}$, and $c_{-\lambda}(\alpha) =
c(-\lambda,\alpha)$ for all $\alpha \in \Delta_+$.

For each $j=1,\dots,n$, set $\omega_j = \alpha_j - \lambda_j$.  For
each $1 \le i \le n-1$ and
$k_i,k_{i+1} \in \mathbb{N}$ with $k_i + k_{i+1} = k \ge 1$, define maps
$$\phi_i = e_{\omega_i}^{\otimes k} \circ (1^{\otimes k_i} \otimes
\mathcal{Y}_c(e^{\lambda_{i-1}},x)^{\otimes k_{i+1}})$$
$$\psi_i = e_{\omega_{i+1}}^{\otimes k} \circ (1^{\otimes k_{i}}
\otimes \mathcal{Y}_c(e^{\lambda_{i+2}},x)^{\otimes k_{i+1}})$$ In the
case that $i=1$, we take $\phi_1 = e_{\omega_1}^{\otimes k}$ and in
the case that $i=n-1$ we take $\psi_{n-1} = e_{\omega_n}^{\otimes k}$.
\begin{theorem}
For every $k_i, k_{i+1} \in \mathbb{N}$ with $k_i + k_{i+1} = k$ and $k \ge 1$, we have
\begin{eqnarray}\label{phimap}
\phi_i:W(k_i\Lambda_i + k_{i+1}\Lambda_{i+1}) \rightarrow
W(k_i\Lambda_0 + k_{i+1}\Lambda_i)
\end{eqnarray}
and
\begin{eqnarray}\label{psimap}
\psi_i:W(k_i\Lambda_i + k_{i+1}\Lambda_{i+1}) \rightarrow
W(k_{i+1}\Lambda_0 + k_{i}\Lambda_{i+1})
\end{eqnarray}
Moreover, for $r_1,\dots , r_n , s \in \mathbb{Z}$,
\begin{eqnarray}\label{phimap'}
\lefteqn{\phi_i:W(k_i\Lambda_i + k_{i+1}\Lambda_{i+1})'_{r_1,\dots ,
    r_n;s}}\\ && \rightarrow W(k_i\Lambda_0 +
k_{i+1}\Lambda_i)'_{r_1,\dots,r_i + k_i, \dots ,
  r_n;s-r_{i-1}+r_i-r_{i+1} + k_i}
\end{eqnarray}
and
\begin{eqnarray}\label{psimap'}
\lefteqn{\psi_i:W(k_i\Lambda_i + k_{i+1}\Lambda_{i+1})'_{r_1,\dots ,
    r_n;s}}\\ && \rightarrow W(k_{i+1}\Lambda_0 +
k_{i}\Lambda_{i+1})'_{r_1,\dots,r_{i+1} + k_{i+1}, \dots ,
  r_n;s-r_{i}+r_{i+1}-r_{i+2} + k_i},
\end{eqnarray}
where we take $r_0 = r_{n+1} = 0$.
\end{theorem}
{\em Proof:} We prove only (\ref{phimap}) since (\ref{psimap}) follows
analogously. Let $a\cdot v_{k_i\Lambda_i + k_{i+1}\Lambda_{i+1}} \in
W(k_i\Lambda_i + k_{i+1}\Lambda_{i+1})$ for some $a \in
U(\bar{\goth{n}})$. We have that
\begin{eqnarray*}
\lefteqn{\phi_i(a\cdot v_{k_i\Lambda_i + k_{i+1}\Lambda_{i+1}})}
\\ &=& \phi_i(a \cdot (\underbrace{e^{\lambda_i} \otimes \dots \otimes
  e^{\lambda_i}}_{k_i-times}\otimes \underbrace{e^{\lambda_{i+1}}
  \otimes \dots \otimes e^{\lambda_{i+1}}}_{k_{i+1}-times}))\\ &=&
\bigg(e_{\omega_i}^{\otimes k} \circ (1^{\otimes k_i} \otimes
\mathcal{Y}_c(e^{\lambda_{i-1}},x)^{\otimes k_{i+1}})\bigg)(a \cdot
(\underbrace{e^{\lambda_i} \otimes \dots \otimes
  e^{\lambda_i}}_{k_i-times}\otimes \underbrace{e^{\lambda_{i+1}}
  \otimes \dots \otimes
  e^{\lambda_{i+1}}}_{k_{i+1}-times}))\\ &=&e_{\omega_i}^{\otimes k}(a
\cdot (\underbrace{e^{\lambda_i} \otimes \dots \otimes
  e^{\lambda_i}}_{k_i-times}\otimes
\underbrace{e_{\lambda_{i-1}}e^{\lambda_{i+1}} \otimes \dots \otimes
  e_{\lambda_{i-1}}e^{\lambda_{i+1}}}_{k_{i+1}-times}))\\ &=&c_1\tau_{\omega_i,c_{-\omega_i}}(a)
\cdot (\underbrace{e^{\alpha_i} \otimes \dots \otimes
  e^{\alpha_i}}_{k_i-times}\otimes \underbrace{e^{\lambda_{i}} \otimes
  \dots \otimes e^{\lambda_{i}}}_{k_{i+1}-times})\\ &=& c_2
\tau_{\omega_i,c_{-\omega_i}}(a) x_{\alpha_i}(-1)^{k_i}\cdot
(\underbrace{\textbf{1} \otimes \dots \otimes
  \textbf{1}}_{k_i-times}\otimes \underbrace{e^{\lambda_{i}} \otimes
  \dots \otimes e^{\lambda_{i}}}_{k_{i+1}-times})\\ &\in&
W(k_i\Lambda_0 + k_{i+1}\Lambda_i)
\end{eqnarray*}
for some constants $c_1,c_2 \in \mathbb{C}$. The fourth equality follows from
the fact that $\lambda_{i-1} + \lambda_{i+1} + \omega_i = \lambda_i$.
This concludes our proof.\\

Using the presentations (\ref{suggested}), we construct exact sequences
which give the multigraded dimensions of certain principal subspaces
(compare to \cite{C1}).

\begin{theorem}\label{shortexact}
Let $k \ge 1$.
For $k_1, k_{2}, k_{n-1}, k_n \in \mathbb{N}$
such that $k_1 + k_{2} = k_{n-1}+ k_n = k$ and $k_1 > 0$ and $k_n > 0$. Then following sequences are exact:
\begin{eqnarray} \label{seq1}
\lefteqn{0 \longrightarrow W(k_1 \Lambda_1 + k_2 \Lambda_2)
  \stackrel{e^{\otimes k}_{\omega_1}} \longrightarrow} \\ && W(k_1
\Lambda_0 + k_2 \Lambda_1) \stackrel{1^{\otimes{k_1-1}} \otimes{\cal
    Y}_{c}(e^{\lambda_1},x)\otimes 1^{\otimes k_{2}}}
\longrightarrow \nonumber \\ && \hspace{2em} W((k_1-1)
\Lambda_0+(k_2+1) \Lambda_1) \longrightarrow 0 \nonumber
\end{eqnarray}
and
\begin{eqnarray}\label{seq2}
\lefteqn{0 \longrightarrow W(k_{n-1} \Lambda_{n -1} + k_{n}
  \Lambda_{n}) \stackrel{e^{\otimes k}_{\omega_n}} \longrightarrow}
\\ && W(k_n \Lambda_0 + k_{n-1} \Lambda_n)
\stackrel{1^{\otimes{k_n-1}} \otimes{\cal
    Y}_{c}(e^{\lambda_n},x)\otimes 1^{\otimes k_{n-1}}}
\longrightarrow \nonumber \\ && \hspace{2em} W((k_n-1)
\Lambda_0+(k_{n-1}+1) \Lambda_n) \longrightarrow 0 \nonumber
\end{eqnarray}
\end{theorem}

{\em Proof:}  We prove that (\ref{seq1}) is exact. The exactness
of (\ref{seq2}) can be proved analogously.
The fact that $e^{\otimes k}_{\omega_1}$ is injective is clear, since its left 
inverse is $e^{\otimes k}_{-\omega_1}$.
We first show that Im$(e^{\otimes k}_{\omega_1}) \subset \mbox{Ker}(1^{\otimes{k_1-1}}
\otimes{\cal Y}_{c}(e^{\lambda_1},x)\otimes 1^{\otimes
  k_{2}})$. Suppose that $w \in \mbox{Im}(e^{\otimes k}_{\omega_1})$. We have that
$$ (1^{\otimes{k_1-1}} \otimes{\cal Y}_{c}(e^{\lambda_1},x)\otimes
1^{\otimes k_{2}})(w) = vx_{\alpha_1}(-1)^{k_1}\cdot v_{(k_1-1)\Lambda_0
  + (k_{2}+1)\Lambda_{1}} = 0
$$ for some $v \in U(\bar{\goth{n}})$, and so $w \in
\mbox{Ker}(1^{\otimes{k_1-1}} \otimes{\cal
  Y}_{c}(e^{\lambda_1},x)\otimes 1^{\otimes k_{2}})$. Hence
Im$(e^{\otimes k}_{\omega_1}) \subset \mbox{Ker}(1^{\otimes{k_1-1}} \otimes{\cal
  Y}_{c}(e^{\lambda_1},x)\otimes 1^{\otimes k_{2}})$.

We now show that $\mbox{Ker}(1^{\otimes{k_1-1}} \otimes{\cal
  Y}_{c}(e^{\lambda_1},x)\otimes 1^{\otimes k_{2}}) \subset
\mbox{Im}(e^{\otimes k}_{\omega_1})$ by characterizing the elements of each set.  If
$w\in \mbox{Ker}(1^{\otimes{k_1-1}} \otimes{\cal
  Y}_{c}(e^{\lambda_1},x)\otimes 1^{\otimes k_{2}})$, we may write
$w = f_{k_1\Lambda_0 + k_{2}\Lambda_1}(u)$ for some $u \in
U(\bar{\goth{n}})$. We have that
$$ (1^{\otimes{k_1-1}} \otimes{\cal Y}_{c}(e^{\lambda_1},x)\otimes
1^{\otimes k_{2}})(f_{k_1\Lambda_0 + k_{2}\Lambda_{k_1}}(u)) =
0\ \mbox{iff} \ f_{(k_{1}-1)\Lambda_0 + (k_{2}+1)\Lambda_1}(u) = 0
$$ 
and by (\ref{1discrepancy}) we have
$$ f_{(k_{1}-1)\Lambda_0 + (k_{2}+1)\Lambda_1}(u) = 0 \ \mbox{iff}
\ u \in \mbox{Ker}f_{(k_1-1)\Lambda_0 + (k_{2}+1)\Lambda_1}
=\mbox{Ker}f_{k \Lambda_0} + U(\bar{\goth{n}})x_{\alpha_1}(-1)^{k_1}
$$ 
so that
$$ w = f_{k_1\Lambda_0 + k_{2}\Lambda_1}(u) \in
\mbox{Ker}(1^{\otimes{k_1-1}} \otimes{\cal
  Y}_{c}(e^{\lambda_1},x)\otimes 1^{\otimes k_{2}}) \ \mbox{iff} \ u
\in \mbox{Ker}f_{(k_1-1)\Lambda_0 + (k_{2}+1)\Lambda_1}.
$$
On the other hand, if $w \in \mathrm{Im} (e^{\otimes k}_{\omega_1})$, we may write
$$ w = vx_{\alpha_1}(-1)^{k_1} \cdot v_{k_1 \Lambda_0 +
  k_{2}\Lambda_1} = f_{k_1\Lambda_0 +
  k_{2}\Lambda_1}(vx_{\alpha_1}(-1)^{k_1})
$$ for some $v \in U(\bar{\goth{n}})$. We may also write
$$
w = f_{k_1 \Lambda_0 + k_{2}\Lambda_1}(u)
$$
for some $u \in U(\bar{\goth{n}})$.
Putting these together, we have that
$$
(u-vx_{\alpha_1}(-1)^{k_1}) \cdot v_{k_1\Lambda_0 + k_{2}\Lambda_1} = 0
$$
which implies
$$
u-vx_{\alpha_1}(-1)^{k_1} \in \mbox{Ker}f_{k_1\Lambda_0 + k_{2}\Lambda_1}.
$$
We therefore have that 
$$ w = f_{k_1\Lambda_0 + k_{2}\Lambda_1}(u) \in \mbox{Im} (e^{\otimes k}_{\omega_1})
\ \mbox{iff} \ u \in \mbox{Ker}f_{k_1\Lambda_0 + k_{2}\Lambda_1} +
U(\bar{\goth{n}})x_{\alpha_1}(-1)^{k_1}
$$
By (\ref{1discrepancy}) we have that 
$$ \mbox{Ker}f_{(k_1-1)\Lambda_0 + (k_{2}+1)\Lambda_1} = \mbox{Ker}f_{k\Lambda_0} +
U(\bar{\goth{n}})x_{\alpha_1}(-1)^{k_1} \subset \mbox{Ker}f_{k_1\Lambda_0 +
  k_{2}\Lambda_1} + U(\bar{\goth{n}})x_{\alpha_1}(-1)^{k_1},
$$
we have that
\begin{eqnarray*}
\lefteqn{w = f_{k_1\Lambda_0 + k_{2}\Lambda_1}(u) \in
  \mbox{Ker}(1^{\otimes{k_1-1}} \otimes{\cal
    Y}_{c}(e^{\lambda_1},x)\otimes 1^{\otimes
    k_{2}})}\\ &\Leftrightarrow& u \in \mbox{Ker}f_{(k_1-1)\Lambda_0 +
  (k_{2}+1)\Lambda_1}\\ &\Rightarrow& u \in \mbox{Ker}f_{k_1\Lambda_0 +
  k_{2}\Lambda_1} +
U(\bar{\goth{n}})x_{\alpha_1}(-1)^{k_1}\\ &\Leftrightarrow& w =
f_{k_1\Lambda_0 + k_{2}\Lambda_1}(u) \in \mbox{Im}(e^{\otimes k}_{\omega_1}),
\end{eqnarray*}
completing our proof.

Using the conjecture presentations (\ref{suggested}), an almost identical proof gives:

\begin{theorem} \label{theorem1}
Let $k \ge 1$.
For any $i$ with $1 \leq i \leq n-1$ and $k_i, k_{i+1} \in \mathbb{N}$
such that $k_i + k_{i+1} = k$, the sequences:
\begin{eqnarray} \label{sequence1}
\lefteqn{W(k_i \Lambda_i + k_{i+1} \Lambda_{i+1}) \stackrel{\phi_i}
  \longrightarrow} \\ && W(k_i \Lambda_0 +k_{i+1} \Lambda_i)
\stackrel{1^{\otimes{k_i-1}} \otimes{\cal
    Y}_{c}(e^{\lambda_i},x)\otimes 1^{\otimes k_{i+1}}}
\longrightarrow \nonumber \\ && \hspace{2em} W((k_i-1)
\Lambda_0+(k_{i+1}+1) \Lambda_i) \longrightarrow 0, \nonumber
\end{eqnarray}
when $k_i \ge 1$, and
\begin{eqnarray} \label{sequence'}
\lefteqn{W(k_i \Lambda_i + k_{i+1} \Lambda_{i+1}) \stackrel{\psi_i}
  \longrightarrow} \\ && W(k_{i+1} \Lambda_0 + k_i \Lambda_{i+1})
\stackrel{1^{\otimes{k_{i+1}-1}} \otimes{\cal
    Y}_{c}(e^{\lambda_{i+1}},x)\otimes 1^{\otimes k_{i}}}
\longrightarrow \nonumber \\ && \hspace{2em} W((k_{i+1}-1)
\Lambda_0+(k_i+1) \Lambda_{i+1}) \longrightarrow 0, \nonumber
\end{eqnarray}
when $k_{i+1} \ge 1$, are exact.
\end{theorem}

\begin{remark} \em
It is important to note that (\ref{seq1}), (\ref{seq2}), (\ref{sequence1}), and (\ref{sequence'})
 are fundamentally different from the
exact sequences used in \cite{C1} and \cite{CalLM2}. In \cite{C1} and
\cite{CalLM2}, exact sequences are constructed using intertwining
operators among level $k$ standard modules. The sequences
(\ref{seq1}), (\ref{seq2}), (\ref{sequence1}), and (\ref{sequence'})
only require intertwining operators among level $1$ standard modules
and recover the same information about multigraded dimensions, as we
will see below. \em
\end{remark}

\begin{remark}\em
Notice that, in general, the first maps in (\ref{sequence1}) and (\ref{sequence'}) are not
injective, in contrast with the corresponding maps in \cite{CLM1}--\cite{CLM2}, \cite{C1}--\cite{C2},
and \cite{CalLM1}--\cite{CalLM3}. In fact, there are only a few cases
where injectivity holds, namely, in the case of Theorem \ref{shortexact} 
and in the corollary below.\em
\end{remark}

\begin{corollary}\label{CalLMseqs}
For each $i=1,\dots,n$ the following sequences are exact:
\begin{eqnarray} 
0 \longrightarrow W(k \Lambda_i) \stackrel{e^{\otimes k}_{\omega_i}}
\longrightarrow W(k \Lambda_0) \stackrel{1^{\otimes{k-1}} \otimes{\cal
    Y}_{c}(e^{\lambda_i},x)} \longrightarrow W((k-1)
\Lambda_0+\Lambda_i) \longrightarrow 0 \nonumber
\end{eqnarray}
\end{corollary}

\begin{remark} \em
These exact sequences in Corollary \ref{CalLMseqs} are the level $k$
analogues of the exact sequences found in \cite{CalLM3}.  \em
\end{remark}

We now use the exact sequences (\ref{seq1}) and (\ref{seq2}) to obtain
the multigraded dimensions $\chi'_{W(k_1\Lambda_1 +
  k_2\Lambda_2)}(x_1,\dots,x_n,q)$ and $\chi'_{W(k_{n-1}\Lambda_{n-1}
  + k_n \Lambda_n)}(x_1,\dots,x_n,q)$.

\begin{theorem}\label{gdimthrm}
Let $k \ge 1$.
Let $k_1, k_2, k_{n-1}, k_n \in \mathbb{N}$ with $k_1 \ge 1$ and $k_n \ge 1$,
such that $k_1 + k_2 = k$ and $k_{n-1} + k_n = k$. Then
\begin{eqnarray} 
\lefteqn{\chi'_{W(k_1\Lambda_1 + k_2\Lambda_2)}(x_1,\dots,x_n,q) =
  \nonumber }\\ &=&x_1^{-k_1}\chi'_{W((k_1-1)\Lambda_0 +
  (k_2+1)\Lambda_1)}(x_1q^{-1},x_2q, x_3\dots,x_n,q)\\ && -
x_1^{-k_1}\chi'_{W(k_1\Lambda_0 +
  k_2\Lambda_1)}(x_1q^{-1},x_2q,x_3,\dots,x_n,q)\nonumber
\end{eqnarray}
and
\begin{eqnarray}
\lefteqn{\chi'_{W(k_{n-1}\Lambda_{n-1} + k_n
    \Lambda_n)}(x_1,\dots,x_n,q) = \nonumber}
\\ &=&x_n^{-k_n}\chi'_{W((k_n-1)\Lambda_0 +
  (k_{n-1}+1)\Lambda_n)}(x_1,\dots,x_{n-1}q,x_nq^{-1},q)\\ && -
x_n^{-k_n}\chi'_{W(k_n\Lambda_0 +
  k_{n-1}\Lambda_n)}(x_1,\dots,x_{n-1}q,x_nq^{-1},q).\nonumber
\end{eqnarray}
\end{theorem}
{\em Proof:} It is easy to see that the maps used in (\ref{seq1}) and (\ref{seq2})
have the property that:
$$ W(k_1 \Lambda_0 + k_2 \Lambda_1)'_{r_1,\dots, r_n,s}
\stackrel{1^{\otimes{k_1-1}} \otimes{\cal
    Y}_{c}(e^{\lambda_1},x)\otimes 1^{\otimes k_{2}}} \longrightarrow
W((k_1-1) \Lambda_0+(k_2+1) \Lambda_1)'_{r_1,\dots,r_n,s}
$$
and
$$ W(k_n \Lambda_0 + k_{n-1} \Lambda_n)'_{r_1,\dots,r_n,s}
\stackrel{1^{\otimes{k_n-1}} \otimes{\cal
    Y}_{c}(e^{\lambda_n},x)\otimes 1^{\otimes k_{n-1}}}
\longrightarrow W((k_n-1) \Lambda_0+(k_{n-1}+1)
\Lambda_n)'_{r_1,\dots,r_n,s}.
$$ Combining this fact with the exactness of (\ref{seq1}) and
(\ref{seq2}), along with (\ref{phimap'}) and (\ref{psimap'}) give
\begin{eqnarray}
\lefteqn{\chi'_{W(k_1\Lambda_0 + k_2\Lambda_1)}(x_1,\dots,x_n,q) = \nonumber }\\
&&x_1^{k_1}q^{k_1}\chi'_{W(k_1\Lambda_1 + k_2\Lambda_2)}(x_1q,x_2q^{-1}, x_3\dots,x_n,q)\\
&& + 
\chi'_{W((k_1-1)\Lambda_0 + (k_2+1)\Lambda_1)}(x_1,\dots,x_n,q)\nonumber
\end{eqnarray}
and
\begin{eqnarray}
\lefteqn{\chi'_{W(k_n\Lambda_0 + k_{n-1}\Lambda_n)}(x_1,\dots,x_n,q) =
  \nonumber }\\ &&x_n^{k_n}q^{k_n}\chi'_{W(k_{n-1}\Lambda_{n-1} +
  k_n\Lambda_n)}(x_1,\dots,x_{n-1}q^{-1},x_nq,q)\\ && +
\chi'_{W((k_n-1)\Lambda_0 +
  (k_{n-1}+1)\Lambda_n)}(x_1,\dots,x_n,q).\nonumber
\end{eqnarray}
which may be rewritten as
\begin{eqnarray}\label{pre1}
\lefteqn{\chi'_{W(k_1\Lambda_1 + k_2\Lambda_2)}(x_1q,x_2q^{-1},
  x_3\dots,x_n,q) = \nonumber}\\ &=&
x_1^{-k_1}q^{-k_1}\chi'_{W(k_1\Lambda_0 +
  k_2\Lambda_1)}(x_1,\dots,x_n,q) \\ &&
-x_1^{-k_1}q^{-k_1}\chi'_{W((k_1-1)\Lambda_0 +
  (k_2+1)\Lambda_1)}(x_1,\dots,x_n,q) \nonumber
\end{eqnarray}
and
\begin{eqnarray}\label{pre2}
\lefteqn{\chi'_{W(k_{n-1}\Lambda_{n-1} + k_n\Lambda_n)}(x_1,\dots,x_{n-1}q^{-1},x_nq,q) = \nonumber}\\
&=&x_n^{-k_n}q^{-k_n}\chi'_{W(k_n\Lambda_0 + k_{n-1}\Lambda_n)}(x_1,\dots,x_n,q)  \\
&&-x_n^{-k_n}q^{-k_n}\chi'_{W((k_n-1)\Lambda_0 + (k_{n-1}+1)\Lambda_n)}(x_1,\dots,x_n,q).\nonumber
\end{eqnarray}
Making the substitutions $$x_1 \mapsto x_1q^{-1},\  x_2 \mapsto x_2q$$ in (\ref{pre1}) and
$$x_n \mapsto x_nq^{-1}, \  x_{n-1} \mapsto x_{n-1}q$$ in (\ref{pre2})
immediately proves our theorem.\\

We now use Theorem \ref{gdimthrm} to write down explicit expressions
for $\chi'_{W(k_1\Lambda_1 + k_2\Lambda_2)}(x_1,\dots,x_n,q)$ and
$\chi'_{W(k_{n-1}\Lambda_{n-1} + k_n\Lambda_n)}(x_1,\dots,x_n,q)$.  In
\cite{G}, Georgiev obtained:
$$
\chi'_{W(k_0 \Lambda_0 + k_j \Lambda_j)}(x_1,\dots,x_n,q) =$$
$$
=\sum_{}\;
\bigg(\frac{q^{{r_{1}^{(1)}}^{2}+ \ldots + {r_{1}^{(k)}}^{2} +
\sum_{t=1}^{k} r_{1}^{(t)}\delta_{1,j_{t}} }}{(q)_{r_{1}^{(1)} -
r_{1}^{(2)}} \ldots (q)_{r_{1}^{(k - 1)} - r_{1}^{(k)}}
(q)_{r_{1}^{(k)}}}\bigg)
\bigg(\frac{q^{{r_{2}^{(1)}}^{2} + \ldots + {r_{2}^{(k)}}^{2} -
r_{2}^{(1)} r_{1}^{(1)} - \ldots - r_{2}^{(k)}r_{1}^{(k)}+
\sum_{t=1}^{k} r_{2}^{(t)}\delta_{2,j_{t}}  }}{(q)_{r_{2}^{(1)}
- r_{2}^{(2)}} \ldots (q)_{r_{2}^{(k - 1)}
- r_{2}^{(k)}} (q)_{r_{2}^{(k)}}}\bigg)\times$$
$$\times\cdots\times \bigg(\frac{q^{{r_{n}^{(1)}}^{2} + \ldots +
    {r_{n}^{(k)}}^{2} - r_{n}^{(1)} r_{n - 1}^{(1)} - \ldots -
    r_{n}^{(k)}r_{n - 1}^{(k)}+ \sum_{t=1}^{k}
    r_{n}^{(t)}\delta_{n,j_{t}} }}{(q)_{r_{n}^{(1)} - r_{n}^{(2)}}
  \ldots (q)_{r_{n}^{(k - 1)} - r_{n}^{(k)}} (q)_{r_{n}^{(k)}}}\bigg)
x_1^{\sum_{i=1}^k r_1^{(i)}}\cdots x_n^{\sum_{i=1}^n r_n^{(i)}}$$
where the sums are taken over decreasing sequences 
$r_j^{(1)} \ge r_j^{(2)} \ge \dots \ge r_j^{(k)} \ge 0$ for each $j=1,\dots,n$
and $j_t=0$ for $0 \le t \le k_0$ and $j_t = j$ for $k_0 < t \le k$,
$j=1,\dots,n$, 
where $(q)_r=\prod_{i=1}^r(1-q^i)$ and $(q)_0 = 1$.
In particular, we have that
$$
\chi'_{W(k_1 \Lambda_0 + k_2 \Lambda_1)}(x_1,\dots,x_n,q) =$$
$$ =\sum_{\begin{array}{c}{\scriptstyle r_{1}^{(1)} \ge\ldots \ge
      r_{1}^{(k)} \geq 0}\\{\scriptstyle 
      \vdots}\\{\scriptstyle r_{n}^{(1)}\ge \ldots \ge
      r_{n}^{(k)} \geq 0 }\end{array}}\; \bigg(\frac{q^{{r_{1}^{(1)}}^{2}+
    \ldots + {r_{1}^{(k)}}^{2} + \sum_{t=k_1+1}^{k} r_{1}^{(t)}
}}{(q)_{r_{1}^{(1)} - r_{1}^{(2)}} \ldots (q)_{r_{1}^{(k - 1)} -
    r_{1}^{(k)}} (q)_{r_{1}^{(k)}}}\bigg)\bigg( \frac{q^{{r_{2}^{(1)}}^{2} +
    \ldots + {r_{2}^{(k)}}^{2} - r_{2}^{(1)} r_{1}^{(1)} - \ldots -
    r_{2}^{(k)}r_{1}^{(k)} }}{(q)_{r_{2}^{(1)} - r_{2}^{(2)}} \ldots
  (q)_{r_{2}^{(k - 1)} - r_{2}^{(k)}} (q)_{r_{2}^{(k)}}}\bigg)\times$$
$$\times\cdots\times
\bigg(\frac{q^{{r_{n}^{(1)}}^{2} + \ldots + {r_{n}^{(k)}}^{2} -
r_{n}^{(1)} r_{n - 1}^{(1)} - \ldots - r_{n}^{(k)}r_{n - 1}^{(k)}
}}{(q)_{r_{n}^{(1)}
- r_{n}^{(2)}} \ldots (q)_{r_{n}^{(k - 1)}
- r_{n}^{(k)}} (q)_{r_{n}^{(k)}}}\bigg) x_1^{\sum_{i=1}^k r_1^{(i)}}\cdots x_n^{\sum_{i=1}^n r_n^{(i)}}$$
and
$$ \chi'_{W(k_n \Lambda_0 + k_{n-1} \Lambda_n)}(x_1,\dots,x_n,q) =$$
$$
=\sum_{\begin{array}{c}{\scriptstyle
r_{1}^{(1)} \ge\ldots \ge  r_{1}^{(k)} \geq 0}\\{\scriptstyle 
\vdots}\\{\scriptstyle r_{n}^{(1)}\ge \ldots
\ge r_{n}^{(k)} \geq 0 }\end{array}}\;
\bigg(\frac{q^{{r_{1}^{(1)}}^{2}+ \ldots + {r_{1}^{(k)}}^{2} 
}}{(q)_{r_{1}^{(1)} -
r_{1}^{(2)}} \ldots (q)_{r_{1}^{(k - 1)} - r_{1}^{(k)}}
(q)_{r_{1}^{(k)}}}\bigg)
\bigg(\frac{q^{{r_{2}^{(1)}}^{2} + \ldots + {r_{2}^{(k)}}^{2} -
r_{2}^{(1)} r_{1}^{(1)} - \ldots - r_{2}^{(k)}r_{1}^{(k)}+
}}{(q)_{r_{2}^{(1)}
- r_{2}^{(2)}} \ldots (q)_{r_{2}^{(k - 1)}
- r_{2}^{(k)}} (q)_{r_{2}^{(k)}}}\bigg) \times$$
$$\times\cdots\times \bigg(\frac{q^{{r_{n}^{(1)}}^{2} + \ldots +
    {r_{n}^{(k)}}^{2} - r_{n}^{(1)} r_{n - 1}^{(1)} - \ldots -
    r_{n}^{(k)}r_{n - 1}^{(k)}+ \sum_{t=k_n+1}^{k} r_{n}^{(t)}
}}{(q)_{r_{n}^{(1)} - r_{n}^{(2)}} \ldots (q)_{r_{n}^{(k - 1)} -
    r_{n}^{(k)}} (q)_{r_{n}^{(k)}}} \bigg) x_1^{\sum_{i=1}^k r_1^{(i)}}\cdots
x_n^{\sum_{i=1}^n r_n^{(i)}}.$$ Applying 
Theorem \ref{gdimthrm} to these expressions immediately gives:
\begin{corollary}\label{gdim}
In the setting of Theorem \ref{gdimthrm}, we have
$$
\chi'_{W(k_1\Lambda_1 + k_2\Lambda_2)}(x_1,\dots,x_n,q)=
$$
$$
=\sum_{}\;
\bigg(\frac{q^{{r_{1}^{(1)}}^{2}+ \ldots + {r_{1}^{(k)}}^{2} +
\sum_{t=k_1+1}^{k} r_{1}^{(t)} + \sum_{t=1}^k r_2^{(t)} - r_1^{(t)} }
(1-q^{r_1^{(k_1)}})}{(q)_{r_{1}^{(1)} -
r_{1}^{(2)}} \ldots (q)_{r_{1}^{(k - 1)} - r_{1}^{(k)}}
(q)_{r_{1}^{(k)}}}\bigg)
\bigg( \frac{q^{{r_{2}^{(1)}}^{2} + \ldots + {r_{2}^{(k)}}^{2} -
r_{2}^{(1)} r_{1}^{(1)} - \ldots - r_{2}^{(k)}r_{1}^{(k)}
 }}{(q)_{r_{2}^{(1)}
- r_{2}^{(2)}} \ldots (q)_{r_{2}^{(k - 1)}
- r_{2}^{(k)}} (q)_{r_{2}^{(k)}}}\bigg) \times$$
$$\times\cdots\times \bigg(
\frac{q^{{r_{n}^{(1)}}^{2} + \ldots + {r_{n}^{(k)}}^{2} -
r_{n}^{(1)} r_{n - 1}^{(1)} - \ldots - r_{n}^{(k)}r_{n - 1}^{(k)}
}}{(q)_{r_{n}^{(1)}
- r_{n}^{(2)}} \ldots (q)_{r_{n}^{(k - 1)}
- r_{n}^{(k)}} (q)_{r_{n}^{(k)}}}\bigg) x_1^{-k_1 + \sum_{i=1}^k r_1^{(i)}}\cdots x_n^{\sum_{i=1}^n r_n^{(i)}}$$
and
$$
\chi'_{W(k_{n-1}\Lambda_{n-1} + k_{n}\Lambda_n)}(x_1,\dots,x_n,q)=
$$
$$
=\sum_{}\;
\bigg(\frac{q^{{r_{1}^{(1)}}^{2}+ \ldots + {r_{1}^{(k)}}^{2} 
}}{(q)_{r_{1}^{(1)} -
r_{1}^{(2)}} \ldots (q)_{r_{1}^{(k - 1)} - r_{1}^{(k)}}
(q)_{r_{1}^{(k)}}}\bigg)
\bigg(\frac{q^{{r_{2}^{(1)}}^{2} + \ldots + {r_{2}^{(k)}}^{2} -
r_{2}^{(1)} r_{1}^{(1)} - \ldots - r_{2}^{(k)}r_{1}^{(k)}+
}}{(q)_{r_{2}^{(1)}
- r_{2}^{(2)}} \ldots (q)_{r_{2}^{(k - 1)}
- r_{2}^{(k)}} (q)_{r_{2}^{(k)}}}\bigg)\times$$
$$\times\cdots\times
\bigg(\frac{q^{{r_{n}^{(1)}}^{2} + \ldots + {r_{n}^{(k)}}^{2} -
r_{n}^{(1)} r_{n - 1}^{(1)} - \ldots - r_{n}^{(k)}r_{n - 1}^{(k)}+
\sum_{t=k_n+1}^{k} r_{n}^{(t)}
}
}{(q)_{r_{n}^{(1)}
- r_{n}^{(2)}} \ldots (q)_{r_{n}^{(k - 1)}
- r_{n}^{(k)}} (q)_{r_{n}^{(k)}}}\bigg) \times
$$
$$\times q^{\sum_{t=1}^kr_{n-1}^{(t)}-r_{n}^{(t)}}(1-q^{r_n^{(k_n)}})
x_1^{\sum_{i=1}^k r_1^{(i)}}\cdots x_n^{-k_n + \sum_{i=1}^n r_n^{(i)}}$$
where the sums are taken over decreasing sequences 
$r_j^{(1)} \ge r_j^{(2)} \ge \dots \ge r_j^{(k)} \ge 0$ for each $j=1,\dots,n$.
\end{corollary}

\begin{remark}
The expressions in Corollary \ref{gdim} can also be written as follows:
As in \cite{G}, for $s= 1,\dots,k-1$ and $i=1,\dots,n$, set $p_i^{(s)} = r_i^{(s)} - r_i^{(s+1)}$,
and set $p_i^{(k)} = r_i^{(k)}$. Also, let $(A_{lm})_{l,m = 1}^n $ be the Cartan matrix of
$\mathfrak{sl}(n+1)$ and $B^{st} := \mathrm{min}\{s,t\}$, $1 \le s,t \le k$. Then,
$$\chi'_{W(k_1\Lambda_1 + k_2\Lambda_2)}(x_1,\dots,x_n,q)=$$
$$\sum_{\begin{array}{c}{\scriptstyle
p_{1}^{(1)} ,\ldots ,  p_{1}^{(k)} \geq 0}\\{\scriptstyle 
\vdots}\\{\scriptstyle p_{n}^{(1)}, \ldots
, p_{n}^{(k)} \geq 0 }\end{array}}
\frac{q^{\frac{1}{2}\sum_{l,m = 1}^n\sum_{s,t=1}^l A_{lm}B^{st}p_l^{(s)} p_m^{(t)}}}{
\prod_{i=1}^n \prod_{s=1}^k (q)_{{p_i}^{(s)}}}q^{\widetilde{p_1}}
q^{\sum_{t=1}^k p_2^{(t)} + \dots + p_2^{(k)} - p_1^{(t)} - \dots - p_1^{(k)}}\times$$
$$\times 
(1-q^{p_1^{(k_1)} + \dots + p_1^{(k)}})x_1^{-k_1}\prod_{i=1}^n x_i^{\sum_{s=1}^k sp_i^{(s)}}
$$
where $\widetilde{p_1} = p_1^{(k_1+1)} + 2p_1^{(k_1+2)} + \dots + k_2p_1^{(k)}$ and 
$$\chi'_{W(k_{n-1}\Lambda_{n-1} + k_n\Lambda_n)}(x_1,\dots,x_n,q)=$$
$$\sum_{\begin{array}{c}{\scriptstyle
p_{1}^{(1)} ,\ldots ,  p_{1}^{(k)} \geq 0}\\{\scriptstyle 
\vdots}\\{\scriptstyle p_{n}^{(1)}, \ldots
, p_{n}^{(k)} \geq 0 }\end{array}}
\frac{q^{\frac{1}{2}\sum_{l,m = 1}^n\sum_{s,t=1}^l A_{lm}B^{st}p_l^{(s)} p_m^{(t)}}}{
\prod_{i=1}^n \prod_{s=1}^k (q)_{{p_i}^{(s)}}}q^{\widetilde{p_n}}
q^{\sum_{t=1}^k p_{n-1}^{(t)} + \dots + p_{n-1}^{(k)} - p_n^{(t)} - \dots - p_n^{(k)}}
\times$$ $$
\times
(1-q^{p_n^{(k_n)} + \dots + p_n^{(k)}}) x_n^{-k_n}\prod_{i=1}^n x_i^{\sum_{s=1}^k sp_i^{(s)}}
$$
where $\widetilde{p_n} = p_n^{(k_n+1)} + 2p_n^{(k_n+2)} + \dots + k_{n-1}p_1^{(k)}$.
\end{remark}

\begin{remark} \em
Corollary \ref{gdim} above is a natural $\widehat{\goth{sl}(n+1)}$-analogue of
Corollary 4.1 in \cite{C1}. The multigraded dimension for
$\chi'_{W(k_1\Lambda_1 + k_2\Lambda_2)}$ in \cite{C1} can be recovered
from the expression above for $\chi'_{W(k_1\Lambda_1 + k_2\Lambda_2)}$
by taking $n=2$. \em
\end{remark}

\begin{remark} \em
Throughout this work we assume that $n\ge2$ for notational convenience.
 In the case that $n=1$ (that is, when $\mathfrak{g} = \mathfrak{sl}(2)$),
 the above results recover the recursions and 
 multigraded dimensions found in \cite{CLM1}-\cite{CLM2}. 
\em
\end{remark}

\setcounter{equation}{0}
\appendix
\section{Appendix}
Above, we needed the construction of a completion of a certain
universal enveloping algebra.  In this section, working in a natural
generality, we recall a construction in \cite{LW} and use it to
construct a natural completion of the universal enveloping algebra of a certain
type of subalgebra of an affine Lie algebras associated to a finite dimensional
semisimple Lie algebra. We also prove a natural decomposition of this completion, which 
is needed above.

Let $\mathfrak{g}$ be a finite dimensional semisimple Lie algebra. Fix
a Cartan subalgebra $\mathfrak{h} \subset \mathfrak{g}$, a set of
roots $\Delta$, a set of simple roots $\Pi =
\{\alpha_1,...,\alpha_n\}$, a set of positive roots
$\Delta_+$, and a symmetric invariant nondegenerate bilinear form
$\langle \cdot, \cdot \rangle$, normalized so that 
$\langle \alpha, \alpha \rangle = 2$ for long roots $\alpha \in \Delta$.
For each $\alpha \in \Delta_+$, let $x_\alpha \in
\mathfrak{g}$ be a root vector associated to the root $\alpha$. We
have that
\begin{equation}\label{strconst}
[x_\alpha, x_\beta] = C_{\alpha, \beta} x_{\alpha + \beta}
\end{equation}
for some constants $C_{\alpha, \beta} \in \mathbb{C}$.  Let $S \subset
\Delta_+$ be a nonempty set of positive roots such that if $\alpha,
\beta \in S$ and $\alpha + \beta \in \Delta_+$, then $\alpha + \beta
\in S$. Define the nilpotent subalgebra $\goth{n}_S \subset
\mathfrak{g}$ by
$$
\goth{n}_S = \sum_{\alpha \in S} \mathbb{C} x_\alpha.
$$
In the case that $S = \Delta_+$, we write $\goth{n}_S = \goth{n}$.

We have the corresponding untwisted affine Lie algebra given by 
\begin{eqnarray*}
\widehat{\mathfrak{g}}= {\mathfrak{g}} \otimes \mathbb{C}[t, t^{-1}]
\oplus \mathbb{C}c,
\end{eqnarray*}
where $c$ is a non-zero central element and 
\begin{eqnarray*}
[ x \otimes t^m, y \otimes t^p ] = [x, y] \otimes t^{m+p} + m\langle x, y \rangle \delta _{m+p, 0} c
\end{eqnarray*}
for any $x, y \in {\mathfrak{g}}$ and $m, p \in \mathbb{Z}$
and 
$$
\bar{\goth{n}}_S  = \goth{n}_S \otimes \mathbb{C}[t,t^{-1}],
$$ a Lie subalgebra of $\widehat{\mathfrak{g}}$. The Lie algebra
$\bar{\goth{n}}_S$ has the following important subalgebras:
$$
\bar{\goth{n}}_{S-} = \goth{n}_S \otimes t^{-1}\mathbb{C}[t^{-1}]
$$
and
$$
\bar{\goth{n}}_{S+} = \goth{n}_S \otimes \mathbb{C}[t].
$$

Let $U(\bar{\goth{n}}_S)$ be the universal 
enveloping algebra of $\bar{\goth{n}}_S$. Using the Poincare-Birkhoff-Witt theorem,
it is easy to see that $U(\bar{\goth{n}}_S)$ has the decomposition
\begin{equation}\label{sdecomp}
U(\bar{\goth{n}}_S) = U(\bar{\goth{n}}_{S-}) \oplus
U(\bar{\goth{n}}_S)\bar{\goth{n}}_{S+}.
\end{equation}

Let $M(S)$ denote the free monoid on $\mathbb{Z} \times S$. We may
write
$$
M(S) = \cup_{n \ge 0} M(S)_n
$$
where
$$
M(S)_n = \mathbb{Z}^n \times S^n
$$
and composition of elements $\circ$ is given by juxtaposition:
$$ (n_1,\dots,n_k; \gamma_1, \dots \gamma_k)\circ (m_1,\dots, m_l ;
\beta_1, \dots , \beta_l) = (n_1,\dots,n_k,m_1,\dots m_l;
\gamma_1,\dots \gamma_k, \beta_1,\dots \beta_l)
$$
where
$$
(n_1,\dots,n_k; \gamma_1, \dots \gamma_k) \in M(S)_k,
$$
$$
(m_1,\dots, m_l ; \beta_1, \dots , \beta_l) \in M(S)_l,
$$
and
$$
(n_1,\dots,n_k,m_1,\dots m_l; \gamma_1,\dots \gamma_k, \beta_1,\dots \beta_l) \in M(S)_{k+l}.
$$ As in \cite{LW}, define, for $n \ge 0$, a map
\begin{eqnarray} \label{tau}
\tau: \mathbb{Z}^n & \longrightarrow & \mathbb{Z}^n \\ \nonumber
(i_1,\dots,i_n) & \mapsto & (i_1+\dots+i_n, i_2+\dots+i_n,\dots,i_n). \nonumber
\end{eqnarray}
For any $b = (n_1,\dots, n_k; \beta_1 , \dots ,\beta_k) \in M(S)_k$
and $i \in \mathbb{Z}$, we write
$$
b \le i \ \mbox{if} \ \tau(n_1,\dots,n_k) \le (i,\dots,i).
$$
In other words, we have
$$
n_1+\dots+n_k \le i,
$$
$$
n_2+\dots+n_k \le i,
$$
$$
\vdots
$$
$$
n_k \le i.
$$

The set $\mbox{Map}(M(S),\mathbb{C})$ of all functions
$$
f: M(S) \longrightarrow \mathbb{C}
$$ has the structure of of an algebra given by taking the identity
element to be the function which is $1$ on $M(S)_0$ and $0$ elsewhere,
and by setting
$$
(r\mu)(a) = r(\mu(a)),
$$
$$
(\mu_1 + \mu_2)(a) = \mu_1(a) + \mu_2(a),
$$
and
$$
(\mu_1\mu_2)(a) = \sum_{a = b \circ c} \mu_1(b) \mu_2(c)
$$ for $r \in \mathbb{C}$, $\mu, \mu_1, \mu_2 \in
\mbox{Map}(M(S),\mathbb{C})$, and $a \in M(S)$.  As in \cite{LW}, for
each $\mu \in \mbox{Map}(M(S),\mathbb{C})$ and $i \in \mathbb{Z}$, we define sets
$$
\mbox{Supp}(\mu) = \{a \in M(S) | \mu(a) \neq 0 \}
$$
and
$$ \mbox{Supp}_i(\mu) = \{ a \in M(\Delta_+) | a \le i \} \cap
\mbox{Supp}(\mu).
$$
Note that if $i \le j$ then $$\mbox{Supp}_i(\mu) \subset \mbox{Supp}_j(\mu)$$
and that 
$$
\mathrm{Supp}(\mu) = \cup_{i \in \mathbb{Z}} \mbox{Supp}_i(\mu).
$$
Define $F(S) \subset \mbox{Map}(M(S),\mathbb{C})$ by
$$
F(S) := \{ \mu:M(S) \longrightarrow \mathbb{C} \  | \  \mathrm{Supp}_i(\mu)\  \mathrm{is \ finite \ 
for\  all\  } i \in \mathbb{Z} \}
$$
and $F_0(S) \subset F(S)$ by
$$
F_0(S) := \{ \mu \in F(S) \  | \  \mbox{Supp}(\mu) \ \mbox{is finite} \}.
$$
We have:
\begin{proposition} (\cite{LW}, Proposition 4.2) 
The set $F(S)$ is a subalgebra of $\mathrm{Map}(M(S),\mathbb{C})$, and $F_0(S) \subset F(S)$ is a 
 subalgebra of $F(S)$. Moreover,  $F_0(S)$ is the free algebra on $\mathbb{Z} \times S$.
\end{proposition}

For each $a \in M(S)$, define maps $X(a) \in F_0(S)$ by
$$
X(a)(b) = \delta_{a,b}.
$$
In particular, for $(n;\beta) \in M(S)_1$, write 
$$
X_\beta(n) = X((n;\beta))
$$
and extend this so that for any $a = (n_1,\dots,n_k;\beta_1,\dots,\beta_k) \in M(S)$
$$
X(a) = X_{\beta_1}(n_1)\dots X_{\beta_k}(n_k).
$$
For any $\mu \in \mbox{Map}(M(S),\mathbb{C})$, we may write
$$
\mu = \sum_{a \in \mbox{Supp}(\mu)} \mu(a)X(a).
$$

Consider the ideal $I_S$ of $F_0(S)$ generated by
\begin{equation*}
 [X_\alpha(n),X_\beta(m)] - C_{\alpha,\beta}X_{\alpha + \beta}(m+n)
\end{equation*}
for $\alpha, \beta \in S$ and $m,n \in \mathbb{Z}$, where
$C_{\alpha,\beta}$ are the structure constants given by (\ref{strconst}). 
\begin{proposition}\label{Uofn}
We have
 $U(\bar{\goth{n}}_S) \simeq F_0(S)/I_S$. 
\end{proposition}
{\em Proof:} Let $T(\bar{\goth{n}}_S)$ denote the tensor algebra over
$\bar{\goth{n}}_S$.  Let $\phi$ be the bijection
\begin{eqnarray} 
\phi: \mathbb{Z} \times S & \longrightarrow & \bar{\goth{n}}_S \\ \nonumber
(n,\beta) & \mapsto & x_{\beta}(n). \nonumber
\end{eqnarray}
Since $F_0(S)$ is the free algebra on $\mathbb{Z} \times \Delta_+$ and
$T(\bar{\goth{n}}_S)$ is the free algebra on $\bar{\goth{n}}_S$, we
extend $\phi$ to a map of free algebras
\begin{eqnarray} 
\phi: F_0(S) & \longrightarrow & T(\bar{\goth{n}}_S) \\ \nonumber
X_{\beta_1}(n_1) \dots X_{\beta_k}(n_k) & \mapsto & x_{\beta_1}(n_1) \dots x_{\beta_k}(n_k), \nonumber
\end{eqnarray}
extended linearly to all of $F_0(S)$. The fact that $\phi$ is
an algebra isomorphism is clear. The proposition follows
immediately.\\

We now impose similar natural relations on $F(S)$. Consider the ideal
$\widetilde{I_S}$ of $F(S)$ generated by
\begin{equation*}
 [X_\alpha(n),X_\beta(m)] - C_{\alpha,\beta}X_{\alpha + \beta}(m+n)
\end{equation*}
for $\alpha, \beta \in S$ and $m,n \in \mathbb{Z}$,
where $C_{\alpha,\beta}$ are the structure constants (\ref{strconst}).
\begin{definition}
 Define the {\em completion} of $U(\bar{\goth{n}}_S)$ by:
 
 \begin{equation} \label{completion}
  \widetilde{U(\bar{\goth{n}}_S)} := F(S) / \widetilde{I_S}.
 \end{equation}
Denote by $[\mu]$ the coset of $\mu \in F(S)$ in $\widetilde{U(\bar{\goth{n}}_S)}$.
\end{definition}

We now introduce some important substructures of
$\widetilde{U(\bar{\goth{n}}_S)}$ and prove some useful facts about
these substructures.  Let
$$ M(S)_- = \{ (m_1,\dots,m_k; \beta_1, \dots, \beta_k) \in M(S) | k
\in \mathbb{N}, m_i \le -1 \mathrm{ \ for \ each \ }i =1,\dots,k \}.
$$

Define
$$ \widetilde{U(\bar{\goth{n}}_{S-})} = \{a \in
\widetilde{U(\bar{\goth{n}}_S)}\ | \ a = [\mu] \mathrm{\ for\ some\ }
\mu \in F(S)\ \mathrm{with\ } \mathrm{Supp}(\mu) \subset M(S)_- \}.
$$
\begin{lemma}
We have that
$$
\widetilde{U(\bar{\goth{n}}_{S-})} \simeq U(\bar{\goth{n}}_{S-}).
$$
\end{lemma}
{\em Proof:} Suppose $[\mu] \in \widetilde{U(\bar{\goth{n}}_{S-})}$
for some $\mu$ with $\mathrm{Supp}(\mu) \subset M(S)_-$. We may write
$$
\mu = \sum_{a \in \mathrm{Supp}(\mu)}\mu(a)X(a)
$$
and so
$$
[\mu] = \sum_{a \in \mathrm{Supp}(\mu)}[\mu(a)X(a)].
$$
By definition, $\mathrm{Supp}_{-1}(\mu)$ is finite, so that 
there are finitely many
$$
a = (m_1,\dots,m_n; \beta_1, \dots, \beta_n) \in \mathrm{Supp}(\mu), \ k \in \mathbb{N}
$$
such that
$$
m_1 + \dots + m_k \le -1
$$
$$
m_2 + \dots + m_k \le -1
$$
$$
\vdots
$$
$$
m_k \le -1.
$$ Since each such $m_i \le -1$, $i=1, \dots k$, have have that
$\mathrm{Supp}_n(\mu) = \mathrm{Supp}_{-1}(\mu)$ for all $n \ge 0$.
In particular, we have that
$$ \mathrm{Supp}(\mu) = \cup_{n \in \mathbb{Z}} \mathrm{Supp}_n(\mu) =
\mathrm{Supp}_{-1}(\mu)
$$ and so $\mathrm{Supp}(\mu)$ is finite and $\mu \in F_0(S)$.  By the
proof of Proposition \ref{Uofn}, we have that
$$
\widetilde{U(\bar{\goth{n}}_{S-})} \simeq U(\bar{\goth{n}}_{S-}),
$$
concluding our proof.\\

Let
$$ M(S)_+ = \{ (m_1,\dots,m_k;\beta_1,\dots,\beta_k)\in M(S)\ | \ k \ge 1\  \mathrm{ and\ there\ exists\ } i \le
k\ \mathrm{with}\ m_i + \dots + m_k \ge 0\}
$$

We define
$$
\widetilde{U(\bar{\goth{n}}_S)\bar{\goth{n}}_{S+}} = \{a \in \widetilde{U(\bar{\goth{n}}_S)} \ | \ 
a = [\mu] \mathrm{\ for\  some\ } \mu \in  F(S)\ 
 \mathrm{with\ } \mathrm{Supp}(\mu) \subset  M(S)_+ \}.
$$
\begin{remark} \em
The space $\widetilde{U(\bar{\goth{n}}_S)\bar{\goth{n}}_{S+}}$ is the
collection of all elements of $\widetilde{U(\bar{\goth{n}}_S)}$ which
have at least one representation as an ``infinite sum" of elements of
$U(\bar{\goth{n}}_S)\bar{\goth{n}}_{S+}$. Indeed, any element $X(a)
\in U(\bar{\goth{n}}_S)$ with $a \in M(S)_+$ can be written as
$$
X(a) = X(b)X(c),
$$
where $i \le k$ and
$$
b = (m_1,\dots,m_{i-1}; \beta_1,\dots,\beta_{i-1}),
$$
$$
c = (m_i,\dots,m_{k};\beta_i,\dots,\beta_k),
$$ and $m_i+\dots +m_k \ge 0$. By (\ref{sdecomp}), $X(c) \in
U(\bar{\goth{n}}_S)\goth{\bar{n}}_{S+}$ and $X(b) \in
U(\bar{\goth{n}}_S)$, and so $X(a) \in
U(\bar{\goth{n}}_S)\goth{\bar{n}}_{S+}$. \em
\end{remark}
\begin{proposition}
The space
$\widetilde{U(\bar{\goth{n}}_S)}$ has the decomposition
\begin{equation}\label{compdecomp}
\widetilde{U(\bar{\goth{n}}_S)} = U(\bar{\goth{n}}_{S-}) \oplus
\widetilde{U(\bar{\goth{n}}_S)\bar{\goth{n}}_{S+}}
\end{equation}
\end{proposition}
{\em Proof:} Given any $u \in U(\bar{\goth{n}}_S)$, using
(\ref{sdecomp}) we may write
$$
u = u_1 + u_2
$$ where $u_1 \in U(\bar{\goth{n}}_{S-})$ and $u_2 \in
U(\bar{\goth{n}}_S)\bar{\goth{n}}_{S+}$.  Suppose $[\mu] \in
\widetilde{U(\bar{\goth{n}}_S)}$ for some $\mu \in F(S)$.  Writing
$$
\mu = \sum_{a \in \mathrm{Supp}(\mu)}\mu(a)X(a),
$$
we have
$$
[\mu] = \sum_{a \in \mathrm{Supp}(\mu)} [\mu(a)X(a)]
$$ and each $[\mu(a)X(a)] \in U(\bar{\goth{n}}_S)$. Since $\mu \in
F(S)$, there are only finitely many $a \in \mathrm{Supp}(\mu)$ such
that $a \in \mathrm{Supp}_{-1}(\mu)$, so that, ranging over all $k\in \mathbb{Z}$, there are only finitely
many $a = (m_1,\dots,m_k; \beta_1, \dots, \beta_k)$ with
$$
m_1 + \dots + m_k \le -1
$$
$$
m_2 + \dots + m_k \le -1
$$
$$
\vdots
$$
$$
m_k \le -1.
$$

For these finitely many $a \in \mathrm{Supp}_{-1}(\mu)$, we write
$$
[\mu(a)X(a)] = [\mu_{1,a}] + [\mu_{2,a}]
$$ for some $[\mu_{1,a}] \in U(\bar{\goth{n}}_{S-})$ and $[\mu_{2,a}]
\in U(\bar{\goth{n}})\bar{\goth{n}}_{S+}$.  By definition of
$U(\bar{\goth{n}}_{S-})$, we have that
$$
\sum_{a \in \mathrm{Supp}_{-1}(\mu)}[\mu_{1,a}] \in U(\bar{\goth{n}}_{S-})
$$ 
since the sum is finite, and
$$ \sum_{a \in \mathrm{Supp}_{-1}(\mu)}[\mu_{2,a}] + \sum_{a \in
  \mathrm{Supp}(\mu) \setminus \mathrm{Supp}_{-1}(\mu)}[\mu(a)X(a)]
\in \widetilde{U(\bar{\goth{n}})\bar{\goth{n}}_{S+}},
$$ since $$\mathrm{Supp}(\mu) \setminus \mathrm{Supp}_{-1}(\mu)
\subset M(S)_+.$$ This shows $[\mu] \in U(\bar{\goth{n}}_{S-}) +
\widetilde{U(\bar{\goth{n}}_S)\bar{\goth{n}}_{S+}}$.  The fact that
$U(\bar{\goth{n}}_{S-}) \cap
\widetilde{U(\bar{\goth{n}}_S)\bar{\goth{n}}_{S+}} = 0$ follows from
the fact that $U(\bar{\goth{n}}_{S-}) \cap
U(\bar{\goth{n}}_S)\bar{\goth{n}}_{S+} = 0$, proving our proposition.

\vspace{.4in}
  
\noindent {\small \sc Department of Mathematics, Rutgers University,
Piscataway, NJ 08854}\\ 
{\em E--mail address}: sadowski@math.rutgers.edu

\end{document}